\documentclass[10pt,leqno,twoside]{article}

\usepackage{graphicx}
\usepackage{color}
\usepackage{comment}

\usepackage[latin2]{inputenc}
\usepackage{amsmath}
\usepackage{amsfonts}
\usepackage{amssymb}
\usepackage{makeidx}
\usepackage{array}
\usepackage{subcaption}

\newtheorem{statement}{Statement}
\newtheorem{theorem}{Theorem}

\newtheorem{lemma}{Lemma}

\newtheorem{remark}[theorem]{Remark}

\newfont{\eurorm}{eurm10 scaled 1100}
\newfont{\eurorms}{eurm10 scaled 800}

\begin{document}

\bigskip
\title{Jordan-Type Inequalities and Stratification}

\noindent
Axioms 13(4), 1--25 (2024), 262\\
https://doi.org/10.3390/axioms13040262

{\let\newpage\relax\maketitle}

\begin{center}
{\em
Milo\v s Mi\' covi\' c${}^{\;\mbox{\scriptsize $\ast$}}$, 
Branko Male\v sevi\' c
\footnote{${\!\!\!\!\!\!\!\!\!\!\!\!}^{\;\mbox{\scriptsize $\ast$}}$correspoding author}
}
\end{center}

\begin{center}
{
\footnotesize
{\it
School of Electrical Engineering, University of Belgrade,                                  \\[-0.00 ex]
Bulevar kralja Aleksandra 73, 11000 Belgrade, Serbia}                                      \\[+0.75 ex]
}
\end{center}

\medskip
\noindent
{\small \textbf{Abstract.}
{\small
$\!\!\!$ 
In this paper, two double Jordan-type inequalities are introduced that generalize some previously established inequalities. As a result, some new upper and lower bounds and approximations of the sinc function are obtained. This extension of Jordan's inequality is enabled by considering the corresponding inequalities through the concept of stratified families of functions. Based on this approach, some optimal approximations of the sinc function are derived by determining the corresponding minimax approximants.
}


\medskip
\noindent
{\small \tt MSC 2020: Primary 41A44; Secondary 26D05}

\section{Introduction}

The function:
$$
\mbox{\rm sinc}\,x 
=
\left\{
\begin{array}{ccc}
 \dfrac{\sin x}{x}, & x \neq 0  \\[1.7 ex]
 1, & x = 0
\end{array}
\right.
$$
has numerous applications in mathematics. The basic approximation of the $\mbox{\rm sinc}\,x$ function is given by the well-known Jordan's inequality:

\begin{theorem}
\label{Theorem_1_0}
{\em \cite{Mitrinovic_1970}}
For $x \!\in\! \left(0, \dfrac{\pi}{2}\right]$, it holds$:$
\begin{equation}
\label{equation_1_0}
\dfrac{2}{\pi}
\leq
\dfrac{\sin x}{x}
<
1 \,.
\end{equation}
\end{theorem}

\medskip\noindent
Since then, many authors have worked on extensions and improvements of Jordan's inequality \cite{Qi_Guo_1993}-\cite{Li_Li_2008}, \cite{Qi_Niu_Guo_1996}-\cite{Popa_2020}. In \cite{Qi_Niu_Guo_1996}, F. Qi, D.-W. Niu and B.-N. Guo did the elaborate research, summarizing previously discovered improvements and applications of Jordan's inequality, along with related problems. Motivated by some of the following results, this paper provides an additional contribution to this topic.

\medskip\noindent
F. Qi and B.-N. Guo, in the paper \cite{Qi_Guo_1993}, provided an enhancement of Jordan's inequality through the following assertion:

\begin{theorem}
\label{Theorem_1_1}
Let $x \!\in\! \left(0, \dfrac{\pi}{2}\right]$. Then, it holds$:$
\begin{equation}
\label{equation_1_1}
\dfrac{2}{\pi} + \dfrac{2}{\pi^2} \left(\pi - 2x \right)
\geq
\dfrac{\sin x}{x}
\geq
\dfrac{2}{\pi} + \dfrac{\pi - 2}{\pi^2} \left(\pi - 2x \right) .
\end{equation}
\end{theorem}

\medskip\noindent
F. Qi then, in the paper \cite{Qi_1996}, provided further improvement of Jordan's inequality through the following assertion:

\begin{theorem}
\label{Theorem_1_2}
Let $x \!\in\! \left(0, \dfrac{\pi}{2}\right]$. Then, it holds$:$
\begin{equation}
\label{equation_1_2}
\dfrac{2}{\pi} + \dfrac{1}{\pi^3} \left(\pi^2 - 4x^2 \right)
\leq
\dfrac{\sin x}{x}
\leq
\dfrac{2}{\pi} + \dfrac{\pi - 2}{\pi^3} \left(\pi^2 - 4x^2 \right) .
\end{equation}
\end{theorem}

\medskip\noindent
In the paper \cite{Deng_1995}, K. Deng contributed to improvements of Jordan's inequality by proving:

\begin{theorem}
\label{Theorem_1_3}
Let $x \!\in\! \left(0, \dfrac{\pi}{2}\right]$. Then, it holds$:$
\begin{equation}
\label{equation_1_3}
\dfrac{2}{\pi} + \dfrac{2}{3\pi^4} \left(\pi^3 - 8x^3 \right)
\leq
\dfrac{\sin x}{x}
\leq
\dfrac{2}{\pi} + \dfrac{\pi - 2}{\pi^4} \left(\pi^3 - 8x^3 \right) .
\end{equation}
\end{theorem}

\medskip\noindent
Based on the inequality (\ref{equation_1_2}), W. D. Jiang and H. Yun provided further extension of Jordan's inequality in their paper \cite{Jiang_Yun_2006} through the following theorem:

\begin{theorem}
\label{Theorem_1_4}
Let $x \!\in\! \left(0, \dfrac{\pi}{2}\right]$. Then, it holds$:$
\begin{equation}
\label{equation_1_4}
\dfrac{2}{\pi} + \dfrac{1}{2\pi^5} \left(\pi^4 - 16x^4 \right)
\leq
\dfrac{\sin x}{x}
\leq
\dfrac{2}{\pi} + \dfrac{\pi - 2}{\pi^5} \left(\pi^4 - 16x^4 \right) .
\end{equation}
\end{theorem}

\medskip\noindent
Shortly afterwards, in the paper \cite{Li_Li_2008}, J.-L. Li and Y.-L. Li provided a more general statement that encompasses the previous inequalities, (\ref{equation_1_1}), (\ref{equation_1_2}), (\ref{equation_1_3}) and (\ref{equation_1_4}), introducing an entire family of inequalities. Namely, the theorem holds:

\begin{theorem}
\label{Theorem_1_5}
Let $x \!\in\! \left(0, \dfrac{\pi}{2}\right]$. Then, it holds$:$
\begin{equation}
\begin{array}{cc}
\label{equation_1_5_1} 
\dfrac{2}{\pi} + \dfrac{2}{\pi^{2}} \left( \pi -  2x \right)
\geq
\dfrac{\sin x}{x}
\geq
\dfrac{2}{\pi} + \dfrac{\pi-2}{\pi^{2}} \left( \pi -  2x \right)
\end{array}
\end{equation}
\begin{equation}
\begin{array}{cc}
\label{equation_1_5_2}
\dfrac{2}{\pi} + \dfrac{2}{n\pi^{n+1}} \left( \pi^{n} \!-\!  \left(2x \right)^{n}\right)
\leq
\dfrac{\sin x}{x}
\leq
\dfrac{2}{\pi} + \dfrac{\pi \!-\! 2}{\pi^{n+1}} \left( \pi^{n} \!-\!  \left(2x \right)^{n}\right) \!\!\! \quad \left(\mbox{for } n \!\in\! N, n \!\geq\! 2 \right).
\end{array}
\end{equation}
\end{theorem}

\medskip\noindent
Inspired by Theorems \ref{Theorem_1_1}, \ref{Theorem_1_2}, \ref{Theorem_1_3}, \ref{Theorem_1_4} and \ref{Theorem_1_5}, in this paper, based on the concept of stratification of corresponding families of functions from the paper \cite{Malesevic_Mihailovic_2021}, we introduce a new extension of Jordan's inequality. Namely, by applying stratification, it is possible to extend the inequality (\ref{equation_1_5_2}) so that the parameter $n$ can be a positive real number. The extension of inequalities for real parameters has recently been the subject of various studies \cite{Malesevic_Micovic_2023}-\cite{Malesevic_Mihailovic_NenezicJovic_Milinkovic_2022}, see also \cite{Chen_Mortici_2023}-\cite{Hung_Li_2016}. Additionally, we provide the best constants for this type of Jordan's inequality, as well as an analysis of upper and lower bounds and minimax approximations of the $\mbox{\rm sinc}\,x$ function based on the inequalities (\ref{equation_1_1}), (\ref{equation_1_2}), (\ref{equation_1_3}), (\ref{equation_1_4}), as well as on the newly obtained inequalities.

\section{Preliminaries}

Recently, in the paper \cite{Malesevic_Mihailovic_2021}, the authors considered families of functions $\varphi_p(x)$, where $x \!\in\! (a,b)  \!\subseteq\! R^{+}$ and $p \!\in\! R^{+}$, which are monotonic with respect to the parameter $p$. In that paper, such families of functions are referred to as stratified families of functions with respect to the parameter $p$.
If, for each $x \!\in\! (a,b)$, it holds:
$$
\left(\forall p_1, p_2 \!\in\! R^{+} \right) \,\,\, p_1 < p_2 \Longleftrightarrow \varphi_{p_1} \left(x \right) < \varphi_{p_2} \left(x \right) ,
$$
then the family of functions $\varphi_p(x)$ is {\em increasingly stratified} with respect to the parameter $p$.
If, for each $x \!\in\! (a,b)$, it holds:
$$
\left(\forall p_1, p_2 \!\in\! R^{+} \right) \,\,\, p_1 < p_2\Longleftrightarrow \varphi_{p_1} \left(x \right) > \varphi_{p_2} \left(x \right) ,
$$
then the family of functions $\varphi_p(x)$ is {\em decreasingly stratified} with respect to the parameter $p$. 

\medskip\noindent
If it is possible to determine a value of the parameter $p=p_0 \!\in\! R^{+}$ for which the infimum of the error:
$$
d_0
=
d(p_0)
=
\sup_{ x \in (a,b)} \left|\varphi_{p_0}(x) \right|
$$
is attained, then the function $\varphi_{p_0}(x)$ is {\em the minimax approximant} of the family of functions $\varphi_{p}(x)$ on the interval $(a, b)$.
Based on the stratifiedness, the parameter value $p = p_0$ is unique.

\medskip\noindent
In this paper, we consider the inequalities (\ref{equation_1_1}), (\ref{equation_1_2}),  (\ref{equation_1_3}),  (\ref{equation_1_4}), (\ref{equation_1_5_1}) and (\ref{equation_1_5_2}) by introducing the corresponding stratified families of functions. When proving inequalities, we will utilize L'H\^ opital's rule for monotonicity, as well as the method for proving MTP (Mixed Trigonometric Polynomial) inequalities described in the paper \cite{Malesevic_Makragic_2016}. 

\medskip\noindent
L'H\^ opital's rule for monotonicity was described by the author I. Pinelis in the paper {\rm $\cite{Pinelis_2002}$}, see also \cite{Estrada_Pavlovic_2017}. In this paper, we use the following formulation:

\begin{lemma} \label{L'Hopital_rule}
$($Monotone form of L'H\^ opital's rule$)$. Let $f$ and $g$ be continuous
functions differentiable on $(a, b)$. Suppose $f(a+) \!=\! g(a+) \!=\! 0$ or
$f(b-) \!=\! g(b-) \!=\! 0$, and assume that $g'(x) \!\neq\! 0$ for all $x \!\in\! (a, b)$. If $f'/g'$
is an increasing (decreasing) function on $(a, b)$, then so is \mbox{$f/g$}.
\end{lemma}

\noindent 
The method to prove inequalities of the form $f(x) > 0$ on the interval $(a, b) \!\subseteq\! R$, where $f(x)$ is an MTP function, as outlined in \cite{Malesevic_Makragic_2016}, is based on determining a downward polynomial approximation $P(x)$ with respect to the observed function $f(x)$. In \cite{Malesevic_Makragic_2016}, the determination of a polynomial $P(x)$ as a polynomial with rational coefficients is considered. If there exists a polynomial  $P(x)$ such that $f(x) > P(x)$ and $P(x) > 0$ on the interval $(a, b)$, then $f(x)>0$ holds on the interval $(a, b)$. The polynomial $P(x) > 0$ is determined as a polynomial with rational coefficients and is examined on the interval $(a, b)$ with rational endpoints. Then, the proof of the inequality $P(x) > 0$ is an algorithmically decidable problem based on Sturm's theorem, see Theorem {\rm 4.2} in \cite{Cutland_1980}. In this paper, the application of Sturm's theorem will not be necessary for proving polynomial inequalities.

\section{Main results}

In this section, several statements are presented and proven, with a special emphasis on the connection between Jordan's inequality and stratification.
Particularly, for each family of functions induced by the aforementioned inequality (\ref{equation_1_5_2}), the best approximations derived from the minimax approximants are identified in Statements \ref{Statement_A} and \ref{Statement_B}.


\begin{lemma}
\label{Lemma_1}
The two-parameter family of functions:
$$
\varphi_{p, q}(x)
=
\dfrac{\sin x}{x} - \dfrac{2}{\pi} - p \left(\pi^q - (2x)^q \right)
$$
is individually decreasingly stratified both with respect to the parameter $p \!\in\! R^{+}$ and with respect to the parameter $q \!\in\! R^{+}$ on the interval $(0, \pi/2)$.
\end{lemma}

\noindent
{\bf Proof.}
For the first derivative of $\varphi_{p, q}(x)$ with respect to $p$, it holds:
$$
\dfrac{\partial {\varphi_{p, q}(x)}}{\partial p}
=
(2x)^q - \pi^q < 0
$$
for $x \!\in\! (0, {\pi}/{2})$ and $q \!\in\! R^{+}$. For the first derivative of $\varphi_{p, q}(x)$ with respect to $q$, it holds:
$$
\dfrac{\partial {\varphi_{p, q}(x)}}{\partial q}
=
p \left( \left(2x \right)^q\ln(2x) - \pi^q\ln(\pi) \right) < 0
$$
for $x \!\in\! (0, {\pi}/{2})$ and $p, q \!\in\! R^{+}$.
\hfill $\Box$

\medskip\noindent
Based on the inequality (\ref{equation_1_5_2}), we introduce the following stratified families of functions in the auxiliary statement:

\begin{lemma}
\label{Lemma_stratification_AB}
Let$:$
$$
A(q)
=
\dfrac{\pi-2}{\pi^{q+1}}
\quad\mbox{and}\quad
B(q)
=
\dfrac{2}{q\pi^{q+1}} \,.
$$
Then, it holds$:$

\medskip
\noindent
{\boldmath $(i)$}
The family of functions$:$
\begin{equation}
\label{family_A}
\varphi_{A(q), q}(x)
=
\dfrac{\sin x}{x} - \dfrac{2}{\pi} - A(q) \left( \pi^{q} - \left(2x\right)^{q} \right)
\end{equation}
is decreasingly stratified with respect to the parameter $q \!\in\! R^{+}$ on the interval $(0, \pi/2)$.

\medskip
\noindent
{\boldmath $(ii)$}
The family of functions$:$
\begin{equation}
\label{family_B}
\varphi_{B(q), q}(x)
=
\dfrac{\sin x}{x} - \dfrac{2}{\pi} - B(q) \left( \pi^{q} - \left(2x\right)^{q} \right)
\end{equation}
is increasingly stratified with respect to the parameter $q \!\in\! R^{+}$ on the interval $(0, \pi/2)$.
\end{lemma}

\noindent
{\bf Proof.}
{\boldmath $(i)$} Since $A(q) = \dfrac{\pi-2}{\pi^{q+1}}$, we obtain the one-parameter family of functions: 
\begin{equation}
\label{equation_A}
\varphi_{A(q), q}(x)
=
\dfrac{\sin x}{x} - 1 +  \left( \dfrac{2x}{\pi}\right) ^{q} \left(1-\dfrac{2}{\pi} \right) .
\end{equation}
The first derivative of $\varphi_{A(q), q}(x)$ with respect to $q$ is:
$$
\dfrac{\partial {\varphi_{A(q), q}(x)}}{\partial q}
=
\left(1-\dfrac{2}{\pi} \right) \left( \dfrac{2x}{\pi}\right) ^{q} \ln\,\dfrac{2x}{\pi} \,.
$$
It is evident that:
$$
\dfrac{\partial {\varphi_{A(q), q}(x)}}{\partial q}  < 0
$$
on the interval $(0, \pi/2)$ for $q \!\in\! R^{+}$, which concludes the proof.

\medskip
\noindent
{\boldmath $(ii)$} Since $B(q) = \dfrac{2}{q\pi^{q+1}}$, we obtain the one-parameter family of
functions:
\begin{equation}
\label{equation_B}
\varphi_{B(q), q}(x)
=
\dfrac{\sin x}{x} - \dfrac{2}{\pi} - \dfrac{2}{q\pi} + \dfrac{2^{q+1}x^q}{q\pi^{q+1}} \,.
\end{equation}
The first derivative of $\varphi_{B(q), q}(x)$ with respect to parameter $q$ is:
$$
\begin{array}{rcl}
\dfrac{\partial {\varphi_{B(q), q}(x)}}{\partial q}
\!\!&\!\!=\!\!&\!\!
\dfrac { 2 } { q^2\pi } + \dfrac {2^{q+1} x^{q} \left( q\ln\,2 + q\ln\,x - q\ln\,\pi - 1 \right)} { q^2\pi^{q+1} } \\[2.5 ex]
\!\!&\!\!=\!\!&\!\!
\dfrac{2}{ q^2\pi } \left(\dfrac{2x}{\pi}\right)^{q} \left( \ln\,\left(\dfrac{2x}{\pi}\right)^{q} + \left(\dfrac{\pi}{2x}\right)^{q} - 1 \right) .
\end{array}
$$
Let $t = \left(\dfrac{2x}{\pi}\right)^{q} $. We now form the function:
$$
g(t) = \ln\left(t \right) + \dfrac{1}{t} - 1 : (0,1) \rightarrow R \,.
$$
Since $\dfrac{d\,g(t)}{d\,t}  = \dfrac{1}{t} -  \dfrac{1}{t^2} < 0$ for $t \!\in\! (0, 1)$, the function $g(t)$ is decreasing on the interval $(0, 1)$. Considering that $g(t)$ is a decreasing function and that $g(1) = 0$,
we conclude that:
$$
g(t) > 0
$$
for $t \!\in\! (0, 1)$.
Thus, it follows:
$$
\dfrac{\partial {\varphi_{B(q), q}(x)}}{\partial q} > 0
$$
on the interval $(0, \pi/2)$ because $g(t)>0$ on $(0, 1)$. This finishes the proof.
\hfill $\Box$


\begin{statement}
\label{Statement_A}
Let$:$
$$
q_1 = \dfrac{2}{\pi-2} = 1.75193 \ldots
\quad\mbox{and}\quad
q_2 = 2 \,.
$$
Then, it holds$:$

\medskip
\noindent
{\boldmath $(i)$}
If $q \!\in\! \left( 0, q_1 \right]$, then the lower bounds of the function $\dfrac{\sin x}{x}$ are given by$:$
$$
x \!\in\! \left(0, \dfrac{\pi}{2}\right)
\;\Longrightarrow\;
\dfrac{\sin x}{x}
>
\dfrac{2}{\pi} + A(q_1) \left( \pi^{q_1} - \left(2x\right)^{q_1} \right)
\geq
\dfrac{2}{\pi} + A(q) \left( \pi^{q} - \left(2x\right)^{q} \right)
$$
and the constant $q_1$ is the best possible.

\medskip
\noindent
{\boldmath $(ii)$}
If $q\!\in\!(q_1, q_2)$, then the equality$:$
\begin{equation*}
\varphi_{A(q), q}(x)
=
\dfrac{\sin x}{x} - \dfrac{2}{\pi} - A(q) \left( \pi^{q} - \left(2x\right)^{q} \right)
=
0
\end{equation*}
has a unique solution $x_0^{(q)}$ and it holds$:$
$$
x\!\in\!\left( 0,x_{0}^{(q)} \right)
\;\Longrightarrow\;
\dfrac{\sin x}{x}
>
\dfrac{2}{\pi} + A(q) \left( \pi^{q} - \left(2x\right)^{q} \right)
$$
and
$$
x\!\in\!\left( x_{0}^{(q)}, \dfrac{\pi}{2}\right)
\;\Longrightarrow\;
\dfrac{\sin x}{x}
<
\dfrac{2}{\pi} + A(q) \left( \pi^{q} - \left(2x\right)^{q} \right) .
$$

\medskip
\noindent
{\boldmath $(iii)$}
If $q \!\in\! \left[ q_2, +\infty \right)$, then the upper bounds of the function $\dfrac{\sin x}{x}$ are given by$:$
$$
x \!\in\! \left(0, \dfrac{\pi}{2}\right)
\;\Longrightarrow\;
\dfrac{\sin x}{x}
<
\dfrac{2}{\pi} + A(q_2) \left( \pi^{q_2} - \left(2x\right)^{q_2} \right)
\leq
\dfrac{2}{\pi} + A(q) \left( \pi^{q} - \left(2x\right)^{q} \right)
$$
and the constant $q_2$ is the best possible.

\medskip
\noindent
{\boldmath $(iv)$}
Each function from the family $\varphi_{A(q), q}(x)$, for $q \!\in\! (q_1, q_2)$, has exactly one maximum and exactly one minimum  at certain points $m_1^{(q)}, m_2^{(q)} \!\in\! (0, \pi/2)$ respectively on the interval $(0, \pi/2)$.
Additionally, it holds $m_1^{(q)} < m_2^{(q)}$. The function $\varphi_{A(q), q}(x)$, for $q=q_1$, has exactly one maximum on $(0, \pi/2)$, and for $q=q_2$ has exactly one minimum  on $(0, \pi/2)$.

\medskip
\noindent
{\boldmath $(v)$}
The equality$:$
$$
\left| \varphi_{A(q), q} \left( m_1^{(q)} \right) \right| = \left| \varphi_{A(q), q} \left( m_2^{(q)} \right) \right|
$$
has the solution $q = q_0$, for the parameter $q \!\in\! (q_1, q_2)$, numerically determined as$:$
$$
q_0
=
1.84823 \ldots \,.
$$
For value$:$
$$
d_{0}
=
\left| \varphi_{A(q_0), q_0} \left( m_1^{(q_0)} \right) \right| = \left| \varphi_{A(q_0), q_0} \left( m_2^{(q_0)} \right) \right|
=
0.0026604 \ldots \,,
$$
it holds$:$
$$
d_0 = \inf\limits_{q \in (0,\infty)} \sup\limits_{x \in (0,\pi/2)}{\left| \varphi_{A(q), q}(x) \right|} \,.
$$
Hence, the minimax approximant of the family of functions $\varphi_{A(q), q}(x)$ is$:$
$$
\varphi_{A(q_0), q_0}(x)
=
\dfrac{\sin x}{x} - \dfrac{2}{\pi} - A(q_0) \left( \pi^{q_0} - \left(2x\right)^{q_0} \right) ,
$$
which determines the corresponding $($minimax$)$ approximation$:$
\begin{equation}
\label{minimax_approx_A}
\dfrac{\sin x}{x}
\approx
\dfrac{2}{\pi} + 0.043803 \ldots \left( \pi^{1.84823 \ldots} - \left(2x\right)^{1.84823 \ldots} \right) .
\end{equation}
\end{statement}

\noindent
{\bf Proof.}
{\boldmath $(i)$}
Let us notice that the assertion is equivalent to ${\varphi_{A(q), q}(x)} > 0$ for $q \leq \dfrac{2}{\pi-2}$ and $x \!\in\! (0, \pi/2)$. Based on $($\ref{equation_A}$)$, it holds: 
\begin{equation}
\label{g_A}
{\varphi_{A(q), q}(x)} = 0
\quad \Longleftrightarrow \quad
q = g(x) = \dfrac{ \ln\,\dfrac{x(\pi-2)}{\pi(x- \sin x)} }{ \ln\,\dfrac{\pi}{2x} } \,.
\end{equation}
We first prove that the function $g(x)$ is monotonic on the interval $(0, \pi/2)$ using L'H\^ opital's rule for monotonicity (Lemma \ref{L'Hopital_rule}).
Let us form the functions $f_1(x) = \ln\,\dfrac{x(\pi-2)}{\pi(x- \sin x)}$ and $f_2(x) = \ln\,\dfrac{\pi}{2x}$ on $(0, \pi/2)$. Note that $f_1(\pi/2-) = 0$ and $f_2(\pi/2-) = 0$. It holds:
$$
\dfrac{f_1'(x)}{f_2'(x)} =\dfrac{-x\cos x+\sin x}{x - \sin x} \,.
$$
We now examine the monotonicity of the function $h(x) = \dfrac{-x\cos x+\sin x}{x - \sin x}$ on the interval $(0, \pi/2)$. The first derivative of the function $h(x)$ is: 
$$
h'(x) =\dfrac{x\cos x+\cos x \sin x + x^2 \sin x - \sin x - x}{ (x - \sin x)^2 } \,.
$$
To examine the sign of the function $h'(x)$, let us examine the sign of the MTP function:
$$
h_1(x) 
= 
x\cos x+\cos x \sin x + x^2 \sin x - \sin x - x
=
x\cos x+ \dfrac{1}{2} \sin 2x + x^2 \sin x - \sin x - x
$$
on the interval $(0, \pi/2)$. \\
We prove that $h_1(x) < 0$ using the method from the paper {\rm \cite{Malesevic_Makragic_2016}}.
If we approximate the functions $\cos x$ and $\sin 2x$ by Maclaurin polynomials of degrees 4 and 9 respectively,
and approximate the function $\sin x$ by Maclaurin polynomial of degree 5 in the addend $x^2 \sin x$
and by Maclaurin polynomial of degree 7 in the addend $- \sin x$, then the function $h_1(x)$ has the upward polynomial approximation:
$$
P_1(x) = \dfrac{2}{2835}x^9 - \dfrac{1}{240}x^7
$$ 
on the interval $(0, \pi/2)$.
It is evident that $P_1(x)<0$ on the interval $(0, \pi/2)$. Thus:
$$
h_1(x) < 0
$$
on the observed interval. From here, we conclude that:
$$
h'(x)<0
$$
on the interval $(0, \pi/2)$. Thus, $h(x) = \dfrac{f_1'(x)}{f_2'(x)}$ is a decreasing function on the interval $(0, \pi/2)$. Furthermore, since $f_1(\pi/2-) = 0$ and $f_2(\pi/2-) = 0$, based on L'H\^ opital's rule for monotonicity, it follows that $g(x) = \dfrac{f_1(x)}{f_2(x)}$ is also a decreasing function on the interval $(0, \pi/2)$. \\
By applying L'H\^ opital's rule, it can be shown that:
$$
\displaystyle{\lim_{x \to \frac{\pi}{2}-}} g(x) = \dfrac{2}{\pi-2} \,.
$$
Considering that $g(x)$ is a decreasing function on the interval $(0, \pi/2)$, we conclude that the function ${\varphi_{A(q), q}(x)}$, for $q = q_1 = \dfrac{2}{\pi-2}$, does not have a root on the observed interval. Since 
$
\displaystyle {\varphi_{A(q_1), q_1}(\pi/4)} = \dfrac{ 2^{\frac{-2}{\pi-2}}(\pi-2) - \pi + 2\sqrt{2} }{\pi} = 0.0082048 \ldots > 0
$, 
we conclude that:
$$
{\varphi_{A(q_1), q_1}(x)}>0
$$
for 
$x \!\in\! (0, \pi/2)$. Additionally, based on the stratification (Lemma \ref{Lemma_stratification_AB}), it holds:
$$
{\varphi_{A(q), q}(x)}>{\varphi_{A(q_1), q_1}(x)}>0
$$
for $q < \dfrac{2}{\pi-2}$ on the interval $(0, \pi/2)$.

\medskip
\noindent
{\boldmath $(ii)$}
It is easily seen that $\displaystyle{\lim_{x \to 0+}} {\varphi_{A(q), q}(x)} = 0$ and $\displaystyle{\lim_{x \to \pi/2-}} {\varphi_{A(q), q}(x)} =0$. In the \mbox{part {\boldmath $(iv)$}} of this proof, it will be shown that each function $\varphi_{A(q), q}(x)$, for $q\!\in\!(q_1, q_2)$, has exactly one maximum and exactly one minimum on the interval $(0, \pi/2)$ respectively. Hence, the stated inequalities follow.

\medskip
\noindent
{\boldmath $(iii)$}
The assertion is equivalent to ${\varphi_{A(q), q}(x)} < 0$ for $q \geq 2$ and $x \!\in\! (0, \pi/2)$.
Continuing from the part {\boldmath $(i)$} of this proof, by multiple applications of L'H\^ opital's rule, it can be shown that:
$$
\displaystyle{\lim_{x \to 0+}} g(x) = 2 \,.
$$
Considering that $g(x)$ is a decreasing function on the interval $(0, \pi/2)$, we conclude that the function ${\varphi_{A(q), q}(x)}$, for $q = q_2= 2$, does not have a root on the observed interval.
Since
$
{\varphi_{A(q_2), q_2}(\pi/4)} = \dfrac{8\sqrt{2}-2-3\pi}{4\pi} = -0.0088386 \ldots < 0
$, 
it holds:
$$
{\varphi_{A(q_2), q_2}(x)}<0
$$
for $x \!\in\! (0, \pi/2)$. Additionally, based on the stratification (Lemma \ref{Lemma_stratification_AB}), it holds:
$$
{\varphi_{A(q), q}(x)}<{\varphi_{A(q_2), q_2}(x)}<0
$$
for $q > 2$ on the interval $(0, \pi/2)$.

\medskip
\noindent
{\boldmath $(iv)$} Let us examine the monotonicity of functions from the family $\varphi_{A(q), q}(x)$ for $q\!\in\!(q_1, q_2)$ on $(0, \pi/2)$. The fourth derivative of $\varphi_{A(q), q}(x)$ with respect to $x$ is:
$$
\dfrac{\partial^4 {\varphi_{A(q), q}(x)}}{\partial x^4}
=
\dfrac{ x^{q+1} f_4(q) + h_4(x) }{x^5} \,,
$$
where
$$
f_4(q)=\pi^{-q-1} 2^q  q(q-1)(q-2)(q-3)(\pi-2)
$$
and
$$
h_4(x)= 4x \left(x^2-6 \right)\cos x + \left(x^4-12x^2+24 \right)\sin x \,.
$$
Moreover, the function $h_4(x)$ is defined at both endpoints of the interval $(0, \pi/2)$, which we will use in the subsequent proof.
The first derivative of the function $h_4(x)$ with respect to $x$ is:
$$
h_4'(x) = x^4\cos x > 0
$$
for $x\!\in\!(0,\pi/2)$. Therefore, the function $h_4(x)$ is increasing on the interval $(0,\pi/2)$. Since $h_4(0) = 0$, it holds that:
$$
h_4(x) > 0
$$
on the interval $(0,\pi/2)$. It is evident that: 
$$
f_4(q) > 0
$$
for $q\!\in\!(q_1, q_2)$. Hence, we have:
\begin{equation}
\label{family_diff_4}
\dfrac{\partial^4 {\varphi_{A(q), q}(x)}}{\partial x^4}
>
0
\end{equation}
on $(0, \pi/2)$ for $q\!\in\!(q_1, q_2)$. Consequently, each function $\dfrac{\partial^3 {\varphi_{A(q), q}(x)}}{\partial x^3}$, for $q\!\in\!(q_1, q_2)$, is increasing on $(0, \pi/2)$. The third derivative of $\varphi_{A(q), q}(x)$ with respect to $x$ is:
$$
\dfrac{\partial^3 {\varphi_{A(q), q}(x)}}{\partial x^3}
=
\dfrac{ x^{q+1}f_3(q) + h_3(x) }{x^4} \,,
$$
where
$$
f_3(q)= \pi^{-q-1} 2^q q (q-1)(q-2)(\pi-2)
\quad\mbox{and}\quad
h_3(x)= \left(-x^3+6x \right)\cos x+ \left(3x^2-6 \right)\sin x \,.
$$
It is evident that $f_3(q) <0$ for $q\!\in\!(q_1, q_2)$.
It holds:
$$
\displaystyle{\lim_{x \to 0+} \! \dfrac{ f_3(q)}{x^{3-q}} \!=\! -\infty} \left(\mbox{for } q\!\in\!(q_1, q_2)\right)
\quad \mbox{and} \quad
\displaystyle{\lim_{x \to 0+} \! \dfrac{h_3(x) }{x^4} \!=\! \lim_{x \to 0+} \! \dfrac{h_3'(x)}{(x^4)'} \!=\! \lim_{x \to 0+} \! \dfrac{x^3\sin x}{4x^3} \!=\! 0 } .
$$
Hence, we have:
\begin{equation}
\label{family_diff_3-1}
\displaystyle{\lim_{x \to 0+} \dfrac{\partial^3 {\varphi_{A(q), q}(x)}}{\partial x^3} = -\infty}
\end{equation}
for $q\!\in\!(q_1, q_2)$. 
It holds:
$$
\displaystyle{\lim_{x \to \frac{\pi}{2}-} \! \dfrac{\partial^3 {\varphi_{A(q), q}(x)}}{\partial x^3} \!=\! \dfrac{(8\pi-16)q^3+(48-24\pi)q^2+(16\pi-32)q+12\pi^2-96}{\pi^4} \!:=\! k_3(q)} \,.
$$
Since $k_3'(q) = \dfrac{8}{\pi^4} (3 q^2- 6 q+2) (\pi-2) > 0$ for $q\!\in\!(q_1, q_2)$, it follows that $k_3(q)$ is an increasing function for $q\!\in\!(q_1, q_2)$. Considering that $k_3(q)$ is an increasing function and that $k_3(q_1) = \dfrac{12\pi^3-48\pi^2-16\pi+160}{\pi^3(\pi-2)^2} = 0.19968 \ldots > 0$, it can be concluded that:
\begin{equation}
\label{family_diff_3-2}
\displaystyle{\lim_{x \to \frac{\pi}{2}-} \dfrac{\partial^3 {\varphi_{A(q), q}(x)}}{\partial x^3} >0}
\end{equation}
for $q\!\in\!(q_1, q_2)$. 
Based on (\ref{family_diff_4}), (\ref{family_diff_3-1}) and (\ref{family_diff_3-2}) each function $\dfrac{\partial^2 {\varphi_{A(q), q}(x)}}{\partial x^2}$, for $q\!\in\!(q_1, q_2)$, has exactly one minimum on $(0, \pi/2)$.
The second derivative of $\varphi_{A(q), q}(x)$ with respect to $x$ is:
$$
\dfrac{\partial^2 {\varphi_{A(q), q}(x)}}{\partial x^2}
=
\dfrac{ x^{q+1}f_2(q) + h_2(x) }{x^3} \,,
$$
where:
$$
f_2(q)= \pi^{-q-1}2^q q (q-1)(\pi-2)
\quad\mbox{and}\quad
h_2(x)= -2x\cos x- \left(x^2-2 \right)\sin x \,.
$$
It is evident that $f_2(q) > 0$ for $q\!\in\!(q_1, q_2)$.
It holds:
$$
\displaystyle{\lim_{x \to 0+} \! \dfrac{ f_2(q)}{x^{2-q}} \!=\! \infty} \left(\mbox{for } q\!\in\!(q_1, q_2)\right)
\quad \!\!\mbox{and}\!\! \quad
\displaystyle{\lim_{x \to 0+} \! \dfrac{h_2(x) }{x^3} \!=\! \lim_{x \to 0+} \! \dfrac{h_2'(x)}{(x^3)'} \!=\! \lim_{x \to 0+} \! \dfrac{-x^2\cos x}{3x^2} \!=\! -\dfrac{1}{3} } .
$$
Hence, we have:
\begin{equation}
\label{family_diff_2-1}
\displaystyle{\lim_{x \to 0+} \dfrac{\partial^2 {\varphi_{A(q), q}(x)}}{\partial x^2} = +\infty}
\end{equation}
for $q\!\in\!(q_1, q_2)$. 
It holds:
$$
\displaystyle{\lim_{x \to \frac{\pi}{2}-} \dfrac{\partial^2 {\varphi_{A(q), q}(x)}}{\partial x^2} = \dfrac{(4\pi-8)q^2+(-4\pi+8)q-2\pi^2+16}{\pi^3} := k_2(q)} \,.
$$
Since $k_2'(q) = \dfrac{4}{\pi^3} (2q - 1) (\pi-2) > 0$ for $q\!\in\!(q_1, q_2)$, it follows that $k_2(q)$ is an increasing function for $q\!\in\!(q_1, q_2)$. Considering that $k_2(q)$ is an increasing function and that $k_2(q_1) = \dfrac{-2\pi^2+4\pi+8}{\pi^2(\pi-2)} = 0.073414 \ldots > 0$, it can be concluded that:
\begin{equation}
\label{family_diff_2-2}
\displaystyle{\lim_{x \to \frac{\pi}{2}-} \dfrac{\partial^2 {\varphi_{A(q), q}(x)}}{\partial x^2} >0}
\end{equation}
for $q\!\in\!(q_1, q_2)$.
We have proven that each function $\dfrac{\partial^2 {\varphi_{A(q), q}(x)}}{\partial x^2}$, for $q\!\in\!(q_1, q_2)$, has exactly one minimum on $(0, \pi/2)$. Therefore, based on (\ref{family_diff_2-1}) and (\ref{family_diff_2-2}), for functions $\dfrac{\partial {\varphi_{A(q), q}(x)}}{\partial x}$, for $q\!\in\!(q_1, q_2)$, there are two possibilities: either they are increasing or they have exactly one maximum and exactly one minimum on $(0, \pi/2)$ respectively. We will prove that:
$$
\displaystyle{\lim_{x \to 0+} \! \dfrac{\partial {\varphi_{A(q), q}(x)}}{\partial x } = 0}
 \mbox{ ,} \quad
\displaystyle{\lim_{x \to \frac{\pi}{2}-} \! \dfrac{\partial {\varphi_{A(q), q}(x)}}{\partial x } > 0}
\quad \mbox{and} \quad
\left( \dfrac{\partial {\varphi_{A(q), q} }}{\partial x} \right) \! (x) \Bigr|_{\substack{x=\frac{\pi}{4}}} < 0
\leqno (\ast)
$$
for $q\!\in\!(q_1, q_2)$, thus, it will be clear that each function $\dfrac{\partial {\varphi_{A(q), q}(x)}}{\partial x}$, for $q\!\in\!(q_1, q_2)$, has exactly one maximum and exactly one minimum on $(0, \pi/2)$ respectively.
The first derivative of $\varphi_{A(q), q}(x)$ with respect to $x$ is:
$$
\dfrac{\partial {\varphi_{A(q), q}(x)}}{\partial x}
=
\dfrac{ x^{q+1}f_1(q) + h_1(x) }{x^2} \,,
$$
where:
$$
f_1(q)= \pi^{-q-1}2^q q (\pi-2)
\quad\mbox{and}\quad
h_1(x)= x\cos x - \sin x \,.
$$
It holds:
$$
\displaystyle{\lim_{x \to 0+} \dfrac{ f_1(q)}{x^{1-q}} = 0} \,\,\, \left(\mbox{for } q\!\in\!(q_1, q_2)\right)
\quad \mbox{and} \quad
\displaystyle{\lim_{x \to 0+} \dfrac{h_1(x) }{x^2} = 0 } \,.
$$
Hence, we have:
\begin{equation}
\label{family_diff_1-2}
\displaystyle{\lim_{x \to 0+} \dfrac{\partial {\varphi_{A(q), q}(x)}}{\partial x } = 0}
\end{equation}
for $q\!\in\!(q_1, q_2)$.
It is easily seen that:
\begin{equation}
\label{family_diff_1-1}
\displaystyle{\lim_{x \to \frac{\pi}{2}-} \dfrac{\partial {\varphi_{A(q), q}(x)}}{\partial x } = \dfrac{2(q(\pi-2)-2)}{\pi^2} > 0}
\end{equation}
for $q\!\in\!(q_1, q_2)$. 
We now examine the sign of the functions $\varphi_{A(q), q}(x)$, for $q\!\in\!(q_1, q_2)$, at the point $x=\pi/4$. It holds that:
\begin{equation*}
\left( \dfrac{\partial {\varphi_{A(q), q} }}{\partial x} \right) (x) \Bigr|_{\substack{x=\frac{\pi}{4}}} = \dfrac{ 2^{-q}q(4\pi-8) +2\sqrt{2} \, (\pi-4)  }{\pi^2} := k_1(q) \,.
\end{equation*}
Since $k_1'(q) = \dfrac{ -4 \, 2^{-q}(\pi-2)(q\ln 2 -1) }{\pi^2} < 0$ for $q\!\in\!(q_1, q_2)$, it follows that $k_1(q)$ is a decreasing function. Considering that $k_1(q)$ is a decreasing function and that $k_1(q_1) = \dfrac{ 2^{ \frac{ 3\pi-8 }{\pi-2} }\pi -  2^{ \frac{ 4\pi-10 }{\pi-2} } + (2\pi^2 - 12\pi + 16)\sqrt{2} }{(\pi-2)\pi^2} = -0.0053418 \ldots < 0$, it can be concluded that:
\begin{equation}
\label{family_diff_1-0}
\left( \dfrac{\partial {\varphi_{A(q), q} }}{\partial x} \right) (x) \Bigr|_{\substack{x=\frac{\pi}{4}}} < 0
\end{equation}
for $q\!\in\!(q_1, q_2)$. Hence, each function $\dfrac{\partial {\varphi_{A(q), q}(x)}}{\partial x}$, for $q\!\in\!(q_1, q_2)$, has exactly one maximum and exactly one minimum on $(0, \pi/2)$ respectively.
Note that $(\ast)$ is a substitution for the conjunction (\ref{family_diff_1-2}), (\ref{family_diff_1-1}), (\ref{family_diff_1-0}).
Additionally, based on the monotonicity of the functions $\dfrac{\partial {\varphi_{A(q), q}(x)}}{\partial x}$, for $q\!\in\!(q_1, q_2)$, and $(\ast)$, we can conclude that each function $\varphi_{A(q), q}(x)$, for $q\!\in\!(q_1, q_2)$, has exactly one maximum and exactly one minimum on $(0, \pi/2)$ respectively. \\
By analyzing the monotonicity of the functions $\dfrac{\partial^4 {\varphi_{A(q), q}(x)}}{\partial x^4}$, $\dfrac{\partial^3 {\varphi_{A(q), q}(x)}}{\partial x^3}$, $\dfrac{\partial^2 {\varphi_{A(q), q}(x)}}{\partial x^2}$, $\dfrac{\partial {\varphi_{A(q), q}(x)}}{\partial x}$ and $\varphi_{A(q), q}(x)$ for $q=q_1$ and for $q=q_2$, in a similar manner, it can be concluded that the function $\varphi_{A(q), q}(x)$, for $q=q_1$, has exactly one maximum on $(0, \pi/2)$, while the function $\varphi_{A(q), q}(x)$, for $q=q_2$, has exactly one minimum on $(0, \pi/2)$.

\medskip
\noindent
{\boldmath $(v)$}
Note that the infimum of the error $d(q) = \sup_{x \in (0,\pi/2)}|\varphi_{A(q), q}(x)|$, for $q\!\in\!(q_1, q_2)$, exists and is attained when:
\begin{equation}
\label{system_A_1}
\left| \varphi_{A(q), q} \left( m_1^{(q)} \right) \right| = \left|  \varphi_{A(q), q} \left( m_2^{(q)} \right) \right| .
\end{equation}
The equation (\ref{system_A_1}) can be numerically solved using the Computer Algebra System \textit{Maple},  yielding in the value of the parameter $q = q_0$ being numerically determined as:
$$
q_0 = 1.84823 \ldots \,,
$$
which determines the minimax approximant $\varphi_{A(q_0), q_0}(x)$ of the family of functions $\varphi_{A(q), q}(x)$. 
\hfill $\Box$

\bigskip\noindent
Figure 1 illustrates the stratified family of functions $\varphi_{A(q), q}$, see $(\ref{family_A})$.
Cases for all values of the parameter $q \!\in\! R^+$ are shown, with a special emphasis on the cases with constants obtained in Statement \ref{Statement_A}.

\medskip\noindent
\begin{figure}[hbt!]
\centering
\includegraphics[width=12cm]{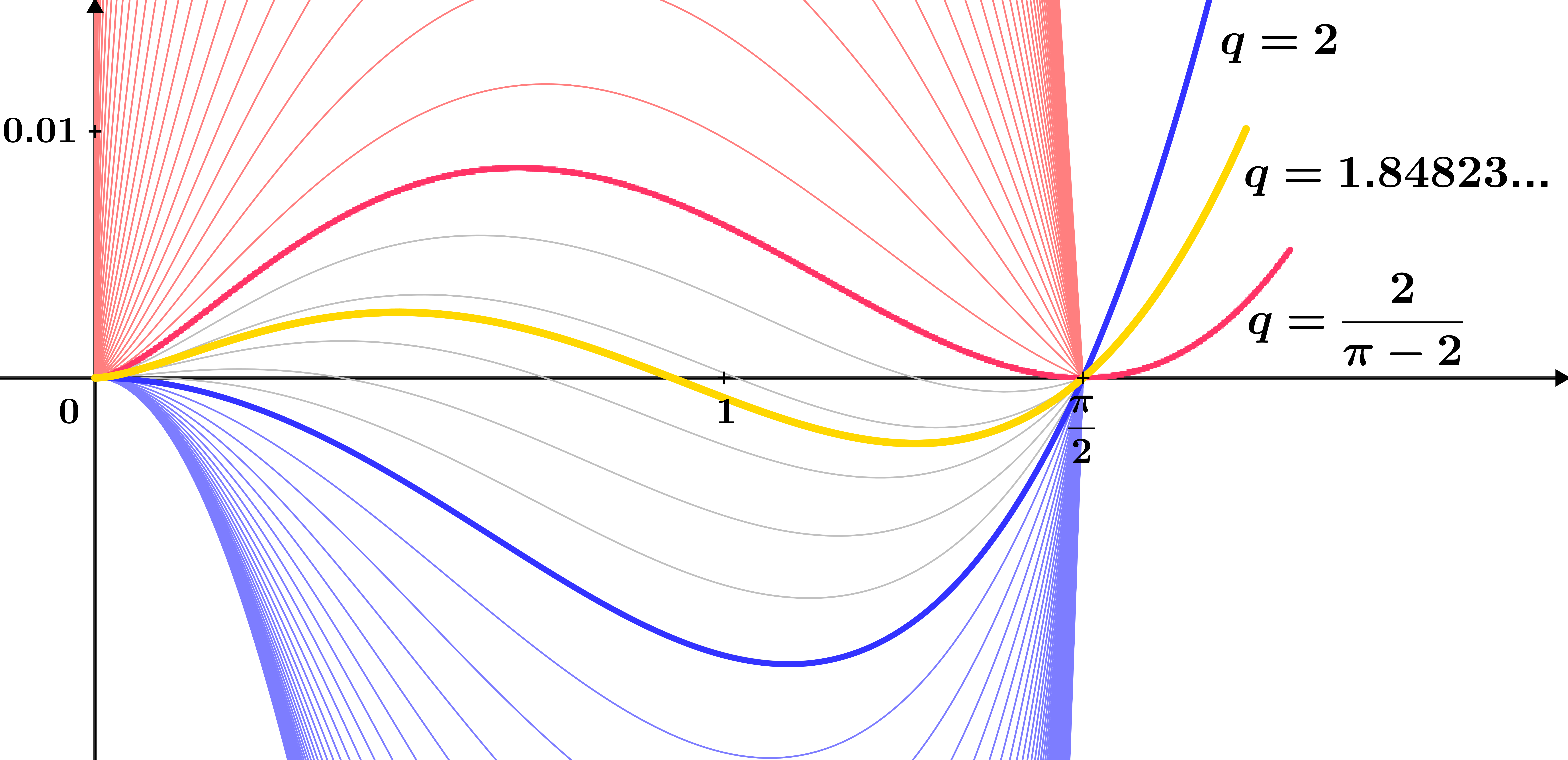}
\caption{Stratified family of functions $\varphi_{A(q), q}$, see $(\ref{family_A})$}
\end{figure}

\break


\begin{statement}
\label{Statement_B}
Let$:$
$$
q_1 = \dfrac{\pi^2}{4}-1 = 1.46740 \ldots
\quad\mbox{and}\quad
q_2 = \dfrac{2}{\pi-2} = 1.75193 \ldots \,.
$$
Then, it holds$:$

\medskip
\noindent
{\boldmath $(i)$}
If $q \!\in\! \left( 0,  q_1 \right]$, then the upper bounds of the function $\dfrac{\sin x}{x}$ are given by$:$
$$
x \!\in\! \left(0, \dfrac{\pi}{2}\right)
\;\Longrightarrow\;
\dfrac{\sin x}{x}
<
\dfrac{2}{\pi} + B(q_1) \left( \pi^{q_1} - \left(2x\right)^{q_1} \right)
\leq
\dfrac{2}{\pi} + B(q) \left( \pi^{q} - \left(2x\right)^{q} \right)
$$
and the constant $q_1$ is the best possible.

\medskip
\noindent
{\boldmath $(ii)$}
If $q\!\in\!(q_1, q_2)$, then the equality$:$
\begin{equation*}
\varphi_{B(q), q}(x)
=
\dfrac{\sin x}{x} - \dfrac{2}{\pi} - B(q) \left( \pi^{q} - \left(2x\right)^{q} \right)
=
0
\end{equation*}
has a unique solution $x_0^{(q)}$ and it holds$:$
$$
x\!\in\!\left( 0,x_{0}^{(q)} \right)
\;\Longrightarrow\;
\dfrac{\sin x}{x}
<
\dfrac{2}{\pi} + B(q) \left( \pi^{q} - \left(2x\right)^{q} \right)
$$
and
$$
x\!\in\!\left( x_{0}^{(q)}, \dfrac{\pi}{2}\right)
\;\Longrightarrow\;
\dfrac{\sin x}{x}
>
\dfrac{2}{\pi} + B(q) \left( \pi^{q} - \left(2x\right)^{q} \right) .
$$

\medskip
\noindent
{\boldmath $(iii)$}
If $q \!\in\! \left[ q_2, +\infty\right)$, then the lower bounds of the function $\dfrac{\sin x}{x}$ are given by$:$
$$
x \!\in\! \left(0, \dfrac{\pi}{2}\right)
\;\Longrightarrow\;
\dfrac{\sin x}{x}
>
\dfrac{2}{\pi} + B(q_2) \left( \pi^{q_2} - \left(2x\right)^{q_2} \right)
\geq
\dfrac{2}{\pi} + B(q) \left( \pi^{q} - \left(2x\right)^{q} \right)
$$
and the constant $q_2$ is the best possible.

\medskip
\noindent
{\boldmath $(iv)$}
Each function from the family $\varphi_{B(q), q}(x)$, for $q \!\in\! (q_1, q_2]$, has exactly one maximum at a point $m^{(q)} \!\in\! (0, \pi/2)$ on the interval $(0, \pi/2)$.

\medskip
\noindent
{\boldmath $(v)$}
The equality$:$
$$
\left| \varphi_{B(q), q} \left( 0+ \right) \right| = \left| \varphi_{B(q), q} \left( m^{(q)} \right) \right|
$$
has the solution $q = q_0$, for the parameter $q \!\in\! (q_1, q_2)$, numerically determined as$:$
$$
q_0
=
1.72287 \ldots \,.
$$
For value$:$
$$
d_{0}
=
\left| \varphi_{B(q_0), q_0} \left( 0+ \right) \right| = \left| \varphi_{B(q_0), q_0} \left( m^{(q_0)} \right) \right|
=
0.0061296 \ldots \,,
$$
it holds$:$
$$
d_0 = \inf\limits_{q \in (0,\infty)} \sup\limits_{x \in (0,\pi/2)}{\left| \varphi_{B(q), q}(x) \right|} .
$$
Hence, the minimax approximant of the family of functions $\varphi_{B(q), q}(x)$ is$:$
$$
\varphi_{B(q_0), q_0}(x)
=
\dfrac{\sin x}{x} - \dfrac{2}{\pi} - B(q_0) \left( \pi^{q_0} - \left(2x\right)^{q_0} \right) ,
$$
which determines the corresponding $($minimax$)$ approximation$:$
\begin{equation}
\label{minimax_approx_B}
\dfrac{\sin x}{x}
\approx
\dfrac{2}{\pi} + 0.051415 \ldots \left( \pi^{1.72287  \ldots} - \left(2x\right)^{1.72287  \ldots} \right) .
\end{equation}
\end{statement}

\noindent
{\bf Proof.}
{\boldmath $(i)$}
Let us notice that the assertion is equivalent to ${\varphi_{B(q), q}(x)} < 0$ for $q \leq \dfrac{\pi^2}{4} - 1$ and $x \!\in\! (0, \pi/2)$. We begin by proving that $\varphi_{B(q), q}(x)$ is monotonic function on the interval $(0, \pi/2)$ for $q = \dfrac{\pi^2}{4} - 1$.
Through elementary transformations, based on $($\ref{equation_B}$)$, it can be shown that the following equivalence holds:
\begin{equation}
\label{g_B}
\begin{array}{cc}
                         &   \dfrac{\partial {\varphi_{B(q), q}(x)}}{\partial x} = \dfrac{x\cos x - \sin x + \left( \dfrac{2x}{\pi} \right)^{q+1}}{x^2} = 0 \\[3.9 ex]
\Longleftrightarrow &
q = g(x) = \dfrac{\ln \dfrac{2x}{ \pi( - x\cos x + \sin x) } }{\ln \dfrac{\pi}{2x}} \,.
\end{array}
\end{equation}
It is necessary to prove that $g(x) \neq \dfrac{\pi^2}{4} - 1$ for every $x \!\in\! (0, \pi/2)$ in order for the function $\varphi_{B(q), q}(x)$ to be monotonic on the interval $(0, \pi/2)$ for $q = \dfrac{\pi^2}{4} - 1$. 
We first prove that the function $g(x)$ is monotonic on the interval $(0, \pi/2)$ by applying L'H\^ opital's rule for monotonicity (Lemma \ref{L'Hopital_rule}).  Let us form the functions $f_1(x) = \ln \dfrac{2x}{ \pi( - x\cos x + \sin x) }$ and $f_2(x) = \ln \dfrac{\pi}{2x}$ on $(0, \pi/2)$. Note that $f_1(\pi/2-) = 0$ and $f_2(\pi/2-) = 0$. It holds:
$$
\dfrac{f_1'(x)}{f_2'(x)} 
=
\dfrac{x\cos x + x^2\sin x - \sin x }{- x\cos x + \sin x} \,.
$$
We now examine the monotonicity of the function $h(x) = \dfrac{x\cos x + x^2\sin x - \sin x }{- x\cos x + \sin x}$ on the interval $(0, \pi/2)$.
The first derivative of the function $h(x)$ is: 
$$
h'(x) =\dfrac{-x (x\cos x \sin x -2\sin^2 x +x^2) }{( - x\cos x + \sin x)^2} \,.
$$
Let us examine the sign of the MTP function:
$$
h_1(x)=x\cos x \sin x -2\sin^2 x +x^2
= \cos 2x + \dfrac{1}{2} x\sin 2x +x^2-1
$$
on the interval $(0, \pi/2)$. If we approximate the functions $\cos 2x$ and $\sin 2x$ by Maclaurin polynomials of degrees 6 and 7  respectively, then the function $h_1(x)$  has the downward polynomial approximation:
$$
P_1(x) = - \dfrac{4}{315}x^8+\dfrac{2}{45}x^6
$$
on the interval $(0, \pi/2)$.
It is evident that $P_1(x)>0$ on the interval $(0, \pi/2)$. Thus:
$$
h_1(x) > 0
$$
on the observed interval. From here, we conclude that:
$$
h'(x)<0
$$
on the observed interval. Thus, $h(x)=\dfrac{f_1'(x)}{f_2'(x)}$ is a decreasing function on the interval $(0, \pi/2)$. Furthermore, since $f_1(\pi/2-) = 0$ and $f_2(\pi/2-) = 0$, based on L'H\^ opital's rule for monotonicity, it follows that $g(x) = \dfrac{f_1(x)}{f_2(x)}$ is also a decreasing function on the interval $(0, \pi/2)$. By applying L'H\^ opital's rule, it can be shown that:
$$
\displaystyle{\lim_{x \to \frac{\pi}{2}-} g(x) = \dfrac{\pi^2}{4} - 1} \,.
$$
Hence, $g(x) > \dfrac{\pi^2}{4} - 1$ on the interval $(0, \pi/2)$. Thus, the function ${\varphi_{B(q), q}(x)}$, for $q = q_1 = \dfrac{\pi^2}{4} - 1$, is monotonic on the interval $(0, \pi/2)$. It holds that
$
\displaystyle{\lim_{x \to 0+} \varphi_{B(q_1), q_1}(x) }= \dfrac{\pi^2-2\pi-4}{\pi^2-4} = -0.070461 \ldots < 0
$
and
$\displaystyle{\lim_{x \to \pi/2-} \varphi_{B(q_1), q_1}(x) = 0 }$. 
Therefore, ${\varphi_{B(q_1), q_1}(x)}$ is an increasing function and negative on $(0, \pi/2)$.
Considering that ${\varphi_{B(q_1), q_1}(x)} < 0$, based on the stratification (Lemma \ref{Lemma_stratification_AB}), it holds:
$$
{\varphi_{B(q), q}(x)} < {\varphi_{B(q_1), q_1}(x)} < 0
$$
for $q < \dfrac{\pi^2}{4} - 1$ on the interval $(0, \pi/2)$.

\smallskip
\noindent
{\boldmath $(ii)$}
Continuing from the previous part of the proof, {\boldmath $(i)$}, by multiple applications of L'H\^ opital's rule, it can be shown that:
$$
\displaystyle{\lim_{x \to 0+}} g(x) = 2 \,.
$$
The function $g(x)$ from (\ref{g_B})  determines the values of the parameter $q$ for which the family of functions $\varphi_{B(q), q}(x)$ have extremes or inflection points on the interval $(0, \pi/2)$. Considering that the function 
$g(x)$ is monotonic on $(0, \pi/2)$  and that $\displaystyle{\lim_{x \to 0+}} g(x) = 2$ and $\displaystyle{\lim_{x \to \pi/2-}} g(x) = \dfrac{\pi^2}{4}-1 = q_1$, every function from the family $\varphi_{B(q), q}(x)$ has either exactly one extremum or exactly one inflection point on the interval $(0, \pi/2)$ for $q\!\in\!\left(\dfrac{\pi^2}{4}-1,2 \right)$, and therefore for $q\!\in\!(q_1, q_2]$, where $q_2 = \dfrac{2}{\pi-2}$, since $q_2<2$.
Let us prove that each function $\varphi_{B(q), q}(x)$, for $q \in (q_1, q_2)$, has exactly one maximum on the interval $(0, \pi/2)$ by proving that all these functions are negative in the right neighborhood of zero and positive and decreasing in the left neighborhood of $\pi/2$.\\
It holds:
$$
\displaystyle{\lim_{x \to 0+}} \varphi_{B(q), q}(x) = \dfrac{(\pi-2)q -2}{\pi q} .
$$ 
Therefore, there exists a right neighborhood of zero such that:
\begin{equation}
\label{taylor_B_0}
\varphi_{B(q), q}(x) < 0
\end{equation}
for $q \in (q_1, q_2)$. The Taylor expansion of the family of functions $\varphi_{B(q), q}(x)$ around $\pi/2$ is:
$$
\varphi_{B(q), q}(x)
=
\dfrac{4q - \pi^2 +4}{\pi^3} \left(x - \dfrac{\pi}{2}\right)^2 + O \left( \left( x - \dfrac{\pi}{2} \right)^3 \right) .
$$
Therefore, there exists a left neighborhood of $\pi/2$ such that:
\begin{equation}
\label{taylor_B_Pi/2}
\varphi_{B(q), q}(x) > 0\quad\mbox{and}\quad
\dfrac{\partial { \varphi_{B(q), q}(x) }}{\partial x} < 0
\end{equation}
for $q \in (q_1, q_2)$. Based on (\ref{taylor_B_0}) and (\ref{taylor_B_Pi/2}) the functions $\varphi_{B(q), q}(x)$, for $q \in (q_1, q_2)$, have exactly one maximum on the interval $(0, \pi/2)$ and the stated inequalities follow.

\medskip
\noindent
{\boldmath $(iii)$}
The assertion is equivalent to ${\varphi_{B(q), q}(x)} > 0$ for $q \geq \dfrac{2}{\pi-2}$ and $x \!\in\! (0, \pi/2)$.
Let us notice that $A(q) = B(q)$ for $q=\dfrac{2}{\pi-2}$, where $A(q) = \dfrac{\pi-2}{\pi^{q+1}}$. In Statement \ref{Statement_A}, it has already been proven that $\varphi_{A(q), q}(x)=\varphi_{B(q), q}(x) > 0$ for $q=q_2=\dfrac{2}{\pi-2}$ on the interval $(0, \pi/2)$. Given that the family of functions $\varphi_{B(q), q}(x)$ is increasingly stratified with respect to the parameter $q$ based on Lemma \ref{Lemma_stratification_AB}, for $q>\dfrac{2}{\pi-2}$, it will also hold that:
$$
\varphi_{B(q), q}(x) > \varphi_{B(q_2), q_2}(x) > 0
$$
on the interval $(0, \pi/2)$.

\medskip
\noindent
{\boldmath $(iv)$} It has been established in the part {\boldmath $(ii)$} of the proof for $q \in (q_1, q_2)$. Similarly, the proof holds for $q=q_2$.

\medskip
\noindent
{\boldmath $(v)$}
Note that the infimum of the error $d(q) = \sup_{x \in (0,\pi/2)}|\varphi_{B(q), q}(x)|$, for $q \!\in\!(q_1, q_2)$, exists and is attained when:
\begin{equation}
\label{system_B_1}
\left| \varphi_{B(q), q} \left( 0+ \right) \right| = \left| \varphi_{B(q), q} \left( m^{(q)} \right) \right| .
\end{equation}
The equation (\ref{system_B_1}) can be numerically solved using the Computer Algebra System \textit{Maple},  yielding in the value of the parameter $q = q_0$ being numerically determined as:
$$
q_0 = 1.72287 \ldots \,,
$$
which determines the minimax approximant $\varphi_{B(q_0), q_0}(x)$ of the family of functions $\varphi_{B(q), q}(x)$. 
\hfill $\Box$

\bigskip\noindent
Figure 2 illustrates the stratified family of functions $\varphi_{B(q), q}$, see $(\ref{family_B})$.
Cases for all values of the parameter $q \!\in\! R^+$ are shown, with a special emphasis on the cases with constants obtained in Statement \ref{Statement_B}.

\medskip\noindent
\begin{figure}[hbt!]
\centering
\includegraphics[width=12cm]{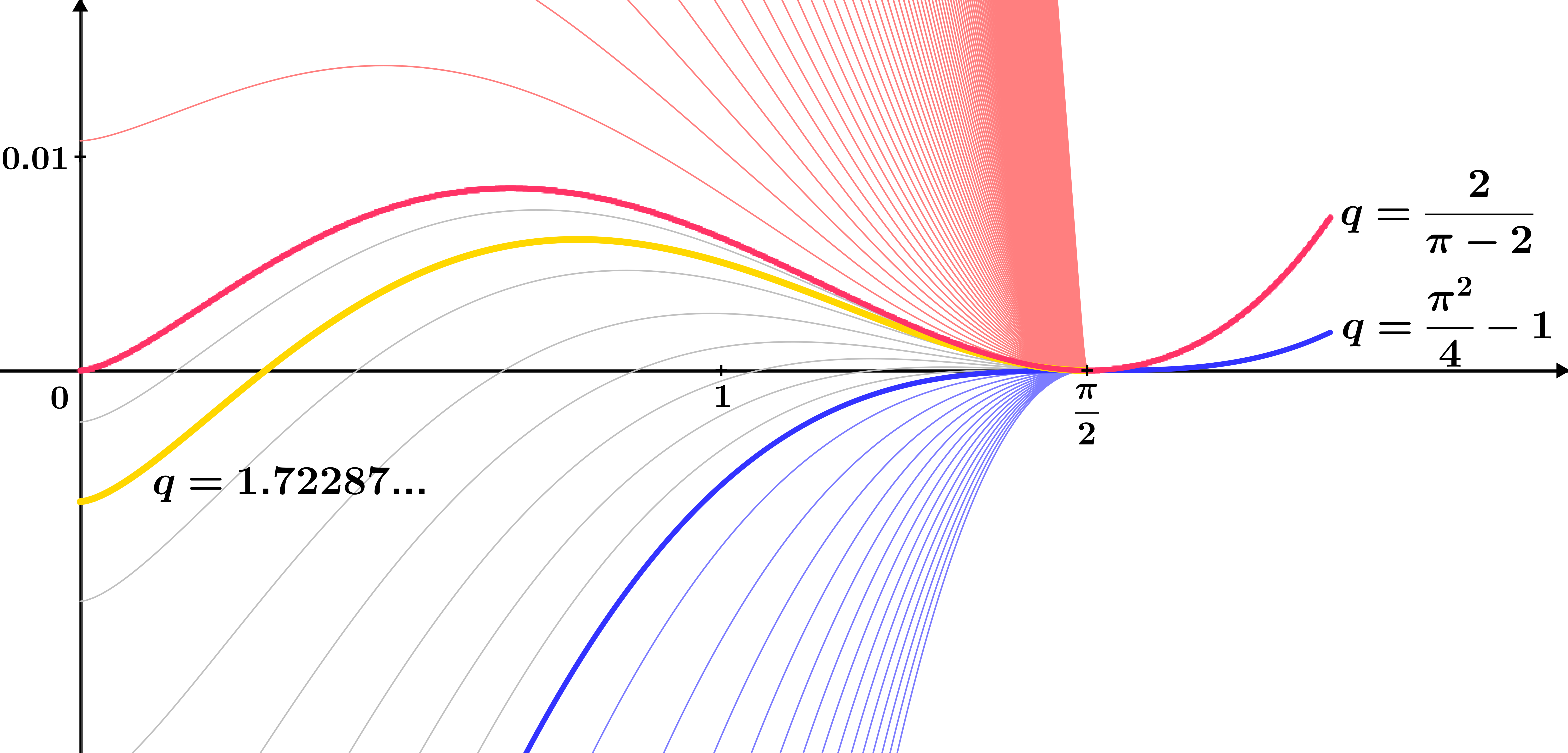}
\caption{Stratified family of functions $\varphi_{B(q), q}$, see $(\ref{family_B})$}
\end{figure}

\break

\bigskip\noindent
In the style of writing Theorem \ref{Theorem_1_5}, based on Statement \ref{Statement_A} and \ref{Statement_B}, we present the following assertion:

\begin{statement}
Let $x \!\in\! \left(0, \dfrac{\pi}{2}\right]$. Then$:$

\medskip
\noindent
{\boldmath $(i)$}
For $q_1 \!\in\! \left( 0,  \dfrac{\pi^2}{4}-1\right] = \left( 0,  1.46740 \ldots \right]$ and $q_2 \!\in\! \left( 0, \dfrac{2}{\pi-2} \right] = \left( 0, 1.75193 \ldots \right]$, it holds$:$
\begin{equation}
\label{equation_main_result_1}
\dfrac{2}{\pi} + \dfrac{2}{q_1\pi^{q_1+1}} \left( \pi^{q_1} - \left(2x\right)^{q_1} \right)
\geq
\dfrac{\sin x}{x}
\geq
\dfrac{2}{\pi} + \dfrac{\pi-2}{\pi^{q_2+1}} \left( \pi^{q_2} - \left(2x\right)^{q_2} \right)  .
\end{equation}

\medskip
\noindent
{\boldmath $(ii)$}
For $q_1 \!\in\! \left[ \dfrac{2}{\pi-2}, +\infty \right) = \left[ 1.75193 \ldots, +\infty \right)$ and $q_2 \!\in\! \left[ 2, +\infty \right)$, it holds$:$
\begin{equation}
\label{equation_main_result_2}
\dfrac{2}{\pi} + \dfrac{2}{q_1\pi^{q_1+1}} \left( \pi^{q_1} - \left(2x\right)^{q_1} \right)
\leq
\dfrac{\sin x}{x}
\leq
\dfrac{2}{\pi} + \dfrac{\pi-2}{\pi^{q_2+1}} \left( \pi^{q_2} - \left(2x\right)^{q_2} \right)  .
\end{equation}
\end{statement}

\begin{remark}
The equalities in {\rm (\ref{equation_main_result_1})} and {\rm (\ref{equation_main_result_2})} clearly hold for $x=\pi/2$.
\end{remark}

\begin{remark}
Note that the inequalities {\rm (\ref{equation_main_result_1})} and {\rm (\ref{equation_main_result_2})} reduce to inequalities {\rm (\ref{equation_1_5_1})} and {\rm (\ref{equation_1_5_2})} respectively when $q_1,q_2\!\in\!N$. 
\end{remark}

\section{Applications}

\medskip\noindent
In this section, we present two applications. The first application is about the improvements and expansions of Theorems \ref{Theorem_1_1},  \ref{Theorem_1_2}, \ref{Theorem_1_3} and \ref{Theorem_1_4}. The second application refers to obtaining some approximations of the $\mbox{\rm sinc}$ function based on some upper and lower bounds of this function and minimax approximants of the corresponding families of functions.

\subsection{Improvements of Theorems \ref{Theorem_1_1}, \ref{Theorem_1_2}, \ref{Theorem_1_3} and \ref{Theorem_1_4}}

\medskip\noindent 
In order to obtain a generalization of all inequalities from Theorems \ref{Theorem_1_2}, \ref{Theorem_1_3}, \ref{Theorem_1_4} and \ref{Theorem_1_5}, for the stratified family of functions $\varphi_{p, q}(x)$ from Lemma \ref{Lemma_1}, we considered the values of the parameter $p=A(q)=\dfrac{\pi-2}{\pi^{q+1}}$ and $p=B(q)=\dfrac{2}{q\pi^{q+1}}$ as functions depending on the parameter $q$. It is possible to consider the family of functions $\varphi_{p, q}(x)$ from Lemma \ref{Lemma_1} by fixing either parameter $p$ or $q$ to some real value. For the cases $q = 1$, $q = 2$, $q = 3$ and $q = 4$, by applying Statement \ref{Statement_A} and \ref{Statement_B}, improvements and extensions of Theorems \ref{Theorem_1_1}, \ref{Theorem_1_2}, \ref{Theorem_1_3}, \ref{Theorem_1_4} respectively can be obtained, as will be shown in the following. 
Particularly, for each family of functions induced by the considered inequalities, 
the best approximations derived from the minimax approximants are identified in Statements \ref{statement_T1}, \ref{statement_T2}, \ref{statement_T3} and \ref{statement_T4}.

\medskip\noindent
In order to improve and extend Theorem \ref{Theorem_1_1}, we consider the family of functions $\varphi_{p, q}(x)$ for the case $q = 1$. The family of functions $\varphi_{p, 1}(x)$ reduces to:
\begin{equation}
\label{family_T1}
\varphi_{p, 1}(x)
=
\dfrac{\sin x}{x} - \dfrac{2}{\pi} - p(\pi - 2x)
\end{equation}
and is decreasingly stratified with respect to the parameter $p \!\in\! R^{+}$ on the interval $(0, \pi/2)$, as proven in Lemma \ref{Lemma_1}.
For this family, the following statement holds:


\begin{statement}
\label{statement_T1}
Let$:$
$$
p_1
=
\dfrac{\pi - 2}{\pi^2}
=
0.11566 \ldots
\quad\mbox{and}\quad
p_2
=
\dfrac{2}{\pi^2}
=
0.20264 \ldots
\,.
$$
Then, it holds$:$

\medskip
\noindent
{\boldmath $(i)$}
If $p\!\in\!(0,p_1]$, then$:$
$$
x \!\in\! \left(0, \dfrac{\pi}{2}\right)
\;\Longrightarrow\;
\dfrac{\sin x}{x}
>
\dfrac{2}{\pi} + p_1(\pi - 2x)
\geq
\dfrac{2}{\pi} + p(\pi - 2x)  \,.
$$

\medskip
\noindent
{\boldmath $(ii)$}
If $p\!\in\!(p_1, p_2)$, then the equality$:$
$$
\varphi_{p, 1}(x)
=
\dfrac{\sin x}{x} - \dfrac{2}{\pi} - p(\pi - 2x)
=
0
$$
has a unique solution $x_0^{(p)}$ and it holds$:$
$$
x\!\in\!\left( 0,x_{0}^{(p)} \right)
\;\Longrightarrow\;
\dfrac{\sin x}{x} < \dfrac{2}{\pi} + p(\pi - 2x)
$$
and
$$
x\!\in\!\left( x_{0}^{(p)}, \dfrac{\pi}{2}\right)
\;\Longrightarrow\;
\dfrac{\sin x}{x} > \dfrac{2}{\pi} + p(\pi - 2x) \,.
$$

\medskip
\noindent
{\boldmath $(iii)$}
If $p\!\in\![p_2,+\infty)$, then$:$
$$
x \!\in\! \left(0, \dfrac{\pi}{2}\right)
\;\Longrightarrow\;
\dfrac{\sin x}{x}
<
\dfrac{2}{\pi} + p_2(\pi - 2x)
\leq
\dfrac{2}{\pi} + p(\pi - 2x) \,.
$$

\medskip
\noindent
{\boldmath $(iv)$}
Each function from the family $\varphi_{p, 1}(x)$, for $p \!\in\! (p_1, p_2]$, has exactly one maximum at a point $m^{(p)} \!\in\! (0, \pi/2)$ on the interval $(0, \pi/2)$.

\medskip
\noindent
{\boldmath $(v)$}
The equality$:$
$$
\left| \varphi_{p, 1} \left(0 + \right)  \right| = \left| \varphi_{p, 1} \left( m^{(p)} \right) \right|
$$
has the solution $p = p_0$, for the parameter $p \!\in\! (p_1, p_2)$, numerically determined as$:$
$$
p_0
=
0.13323 \ldots \,.
$$
For value$:$
$$
d_{0}
=
\left| \varphi_{p_0, 1} \left( 0 + \right) \right| = \left| \varphi_{p_0, 1} \left(m^{(p_0)} \right) \right|
=
0.055187 \ldots \,,
$$
it holds$:$
$$
d_0 = \inf\limits_{p \in (0,\infty)} \sup\limits_{x \in (0,\pi/2)}{\left| \varphi_{p, 1}(x) \right|} \, .
$$
Hence, the minimax approximant of the family of functions $\varphi_{p, 1}(x)$ is$:$
$$
\varphi_{p_0, 1}(x)
=
\dfrac{\sin x}{x} - \dfrac{2}{\pi} - p_0(\pi - 2x) \,,
$$
which determines the corresponding $($minimax$)$ approximation$:$
\begin{equation}
\label{minimax_approx_T1}
\dfrac{\sin x}{x}
\approx
\dfrac{2}{\pi} + 0.13323 \ldots (\pi - 2x) \, .
\end{equation}
\end{statement}

\noindent
{\bf Proof.}
{\boldmath $(i)$} The claim follows directly from Statement \ref{Statement_A} and based on the stratification. Namely, for $q=1$, it holds that $A(q) = \dfrac{\pi-2}{\pi^{q+1}} = p_1$.

\medskip
\noindent
{\boldmath $(ii)$}
Let us examine the monotonicity of functions ${\varphi_{p, 1}(x)}$ for $p \!\in\! (p_1, p_2)$ on the interval $(0, \pi/2)$ in a similar manner as in the proof of Statement \ref{Statement_A}.
The second derivative of $\varphi_{p, 1}(x)$ with respect to $x$ is:
$$
\dfrac{\partial^2 {\varphi_{p, 1}(x)}}{\partial x^2}
=
\dfrac{f(x)}{x^{3}} \,,
$$
where the function $f(x)$ is an MTP function given by:
$$
f(x)
=
-2x\cos x -x^{2}\sin x +2\sin x \,.
$$
Let us note that:
$$
f'(x) = -x^{2}\cos x < 0
$$
on the interval $(0, \pi/2)$. Thus, the function $f(x)$ is decreasing on the observed interval. Considering that $f(x)$ is a decreasing function
on the interval $(0, \pi/2)$ and that $f(0+)=0$, it follows that:
$$
f(x)<0
$$
for $x \!\in\! (0, \pi/2)$. Hence:
\begin{equation}
\label{app1_diff_3}
\dfrac{\partial^2 {\varphi_{p, 1}(x)}}{\partial x^2} < 0
\end{equation}
for $x \!\in\! (0, \pi/2)$.\\
The Taylor expansion of the family of functions $\varphi_{p, 1}(x)$ around zero is:
$$
\varphi_{p, 1}(x)
=
\left(1 -  \dfrac{2}{\pi} - p\pi \right) +2px + O\left(x^2 \right) .
$$
Therefore, there exists a right neighborhood of zero such that:
\begin{equation}
\label{app1_diff_0}
\varphi_{p, 1}(x) < 0
\quad\mbox{and}\quad
\dfrac{\partial {\varphi_{p, 1}(x)}}{\partial x} > 0 
\end{equation}
for $p\!\in\!(p_1, p_2)$.
The Taylor expansion of the family of functions $\varphi_{p, 1}(x)$ around $\pi/2$ is:
$$
\varphi_{p, 1}(x)
=
\left( -\dfrac{4}{\pi^2}+2 p \right) \left(x - \dfrac{\pi}{2}\right) + O \left( \left( x - \dfrac{\pi}{2} \right)^2 \right) .
$$
Therefore, there exists a left neighborhood of $\pi/2$ such that:
\begin{equation}
\label{app1_diff_Pi/2}
\varphi_{p, 1}(x) > 0
\quad\mbox{and}\quad
\dfrac{\partial {\varphi_{p, 1}(x)}}{\partial x} < 0
\end{equation}
for $p\!\in\!(p_1, p_2)$.\\
By analyzing the monotonicity of the functions $\dfrac{\partial^2 {\varphi_{p, 1}(x)}}{\partial x^2}$, $\dfrac{\partial {\varphi_{p, 1}(x)}}{\partial x}$ and $\varphi_{p, 1}(x)$ for $p \!\in\! (p_1, p_2)$ on the interval $(0, \pi/2)$, in a similar manner as in the proof of Statement \ref{Statement_A}, based on (\ref{app1_diff_3}), (\ref{app1_diff_0}) and (\ref{app1_diff_Pi/2}),
it can be concluded that each function 
$\varphi_{p, 1}(x)$, for $p \!\in\! (p_1, p_2)$, has exactly one maximum on the interval $(0, \pi/2)$. From $\displaystyle{\lim_{x \to 0+} \varphi_{p, 1}(x)} < 0$ and $\displaystyle{\lim_{x \to \pi/2-} \varphi_{p, 1}(x)} > 0$, for $p\!\in\!(p_1, p_2)$, the corresponding inequalities follow.

\medskip
\noindent
{\boldmath $(iii)$}
The claim follows directly from Statement \ref{Statement_B} and based on the stratification. Namely, for $q=1$, it holds that $B(q)=\dfrac{2}{q\pi^{q+1}} = p_2$.

\medskip
\noindent
{\boldmath $(iv)$}
It has been proven within the proof {\boldmath $(ii)$}.

\medskip
\noindent
{\boldmath $(v)$}
Note that the infimum of the error $d(p) = \sup_{ x \in (0,\pi/2)}|\varphi_{p, 1}(x)|$, for $p \!\in\! (p_1, p_2)$, exists and is attained when:
\begin{equation}
\label{system_app0_1}
\left| \varphi_{p, 1} \left(0 + \right) \right| = \left| \varphi_{p, 1} \left( m^{(p)} \right) \right| .
\end{equation}
The equation (\ref{system_app0_1}) can be numerically solved using the Computer Algebra System \textit{Maple},  yielding in the value of the parameter $p = p_0$ being numerically determined as:
$$
p_0 = 0.13323\ldots \,,
$$
which determines the minimax approximant $\varphi_{p_0, 1}(x)$ of the family of functions $\varphi_{p, 1}(x)$.
\hfill $\Box$
\noindent

\medskip\noindent
In order to improve and extend Theorem \ref{Theorem_1_2}, we consider the family of functions $\varphi_{p, q}(x)$ for the case $q = 2$. The family of functions $\varphi_{p, 2}(x)$ reduces to:
\begin{equation}
\label{family_T2}
\varphi_{p, 2}(x)
=
\dfrac{\sin x}{x} - \dfrac{2}{\pi} - p \left(\pi^2 - 4x^2 \right)
\end{equation}
and is decreasingly stratified with respect to the parameter $p \!\in\! R^{+}$ on the interval $(0, \pi/2)$, as proven in Lemma \ref{Lemma_1}.
For this family, the following statement holds:


\begin{statement}
\label{statement_T2}
Let$:$
$$
p_1
=
\dfrac{1}{\pi^3}
=
0.032251 \ldots
\quad\mbox{and}\quad
p_2
=
\dfrac{\pi - 2}{\pi^3}
=
0.036818 \ldots \,.
$$
Then, it holds$:$

\medskip
\noindent
{\boldmath $(i)$}
If $p\!\in\!(0,p_1]$, then$:$
$$
x \!\in\! \left(0, \dfrac{\pi}{2}\right)
\;\Longrightarrow\;
\dfrac{\sin x}{x}
>
\dfrac{2}{\pi} + p_1 \left(\pi^2 - 4x^2 \right)
\geq
\dfrac{2}{\pi} + p \left(\pi^2 - 4x^2 \right) .
$$

\medskip
\noindent
{\boldmath $(ii)$}
If $p\!\in\!(p_1, p_2)$, then the equality$:$
$$
\varphi_{p, 2}(x)
=
\dfrac{\sin x}{x} - \dfrac{2}{\pi} - p \left(\pi^2 - 4x^2 \right)
=
0
$$
has a unique solution $x_0^{(p)}$ and it holds$:$
$$
x\!\in\!\left( 0,x_{0}^{(p)} \right)
\;\Longrightarrow\;
\dfrac{\sin x}{x} > \dfrac{2}{\pi} + p \left(\pi^2 - 4x^2 \right)
$$
and
$$
x\!\in\!\left( x_{0}^{(p)}, \dfrac{\pi}{2}\right)
\;\Longrightarrow\;
\dfrac{\sin x}{x} < \dfrac{2}{\pi} + p \left(\pi^2 - 4x^2 \right) .
$$

\medskip
\noindent
{\boldmath $(iii)$}
If $p\!\in\![p_2,+\infty)$, then$:$
$$
x \!\in\! \left(0, \dfrac{\pi}{2}\right)
\;\Longrightarrow\;
\dfrac{\sin x}{x}
<
\dfrac{2}{\pi} + p_2 \left(\pi^2 - 4x^2 \right)
\leq
\dfrac{2}{\pi} + p \left(\pi^2 - 4x^2 \right) .
$$

\medskip
\noindent
{\boldmath $(iv)$}
Each function from the family $\varphi_{p, 2}(x)$, for $p \!\in\! (p_1, p_2]$, has exactly one minimum at a point $m^{(p)} \!\in\! (0, \pi/2)$ on the interval $(0, \pi/2)$.

\medskip
\noindent
{\boldmath $(v)$}
The equality$:$
$$
\left| \varphi_{p, 2} \left(0 + \right) \right| = \left| \varphi_{p, 2} \left( m^{(p)} \right) \right|
$$
has the solution $p = p_0$, for the parameter $p \!\in\! (p_1, p_2)$, numerically determined as$:$
$$
p_0
=
0.036014 \ldots \,.
$$
For value$:$
$$
d_{0}
=
\left| \varphi_{p_0, 2} \left( 0 + \right) \right| = \left| \varphi_{p_0, 2} \left(m^{(p_0)} \right) \right|
=
0.0079283 \ldots \,,
$$
it holds$:$
$$
d_0 = \inf\limits_{p \in (0,\infty)} \sup\limits_{x \in (0,\pi/2)}{\left| \varphi_{p, 2}(x) \right|} \,.
$$
Hence, the minimax approximant of the family of functions $\varphi_{p, 2}(x)$ is$:$
$$
\varphi_{p_0, 2}(x)
=
\dfrac{\sin x}{x} - \dfrac{2}{\pi} - p_0 \left(\pi^2 - 4x^2 \right) ,
$$
which determines the corresponding $($minimax$)$ approximation$:$
\begin{equation}
\label{minimax_approx_T2}
\dfrac{\sin x}{x}
\approx
\dfrac{2}{\pi} + 0.036014 \ldots \left(\pi^2 - 4x^2 \right) .
\end{equation}
\end{statement}

\noindent
{\bf Proof.}
{\boldmath $(i)$} The claim follows directly from Statement \ref{Statement_B} and based on the stratification. Namely, for $q=2$, it holds that $B(q) = \dfrac{2}{q\pi^{q+1}} = p_1$.

\medskip
\noindent
{\boldmath $(ii)$}
Let us examine the monotonicity of functions ${\varphi_{p, 2}(x)}$ for $p \!\in\! (p_1, p_2)$ on the interval $(0, \pi/2)$ in a similar manner as in the proof of Statement \ref{Statement_A}.
 The third derivative of $\varphi_{p, 2}(x)$ with respect to $x$ is:
$$
\dfrac{\partial^3 {\varphi_{p, 2}(x)}}{\partial x^3}
=
\dfrac{f(x)}{x^{4}} \,,
$$
where the function $f(x)$ is an MTP function given by:
$$
f(x)
=
-x^{3}\cos x + 6x\cos x + 3x^{2}\sin x - 6\sin x \,.
$$
Let us note that:
$$
f'(x) = x^{3}\sin x > 0
$$
on the interval $(0, \pi/2)$. Thus, the function $f(x)$ is increasing on the observed interval. Considering that $f(x)$ is an increasing function
on the interval $(0, \pi/2)$ and that $f(0+)=0$, it follows that:
$$
f(x)>0
$$
for $x \!\in\! (0, \pi/2)$. Hence:
\begin{equation}
\label{app_diff_3}
\dfrac{\partial^3 {\varphi_{p, 2}(x)}}{\partial x^3} > 0
\end{equation}
for $x \!\in\! (0, \pi/2)$.\\
The Taylor expansion of the family of functions $\varphi_{p, 2}(x)$ around zero is:
$$
\varphi_{p, 2}(x)
=
\left(1 -  \dfrac{2}{\pi} - p\pi^2 \right) + \left(- \dfrac{1}{6} + 4p \right)x^2 + O\left(x^4 \right) .
$$
Therefore, there exists a right neighborhood of zero such that:
\begin{equation}
\label{app_diff_0}
\varphi_{p, 2}(x) > 0
\mbox{ ,}\quad
\dfrac{\partial {\varphi_{p, 2}(x)}}{\partial x} < 0
\quad\mbox{and}\quad
\dfrac{\partial^2 {\varphi_{p, 2}(x)}}{\partial x^2} < 0
\end{equation}
for $p\!\in\!(p_1, p_2)$.
The Taylor expansion of the family of functions $\varphi_{p, 2}(x)$ around $\pi/2$ is:
$$
\varphi_{p, 2}(x)
=
\left( -\dfrac{4}{\pi^2}+4\pi p \right) \left(x - \dfrac{\pi}{2}\right) + \left( \dfrac{8}{\pi^3} - \dfrac{1}{\pi}+ 4p\right) \left(x - \dfrac{\pi}{2}\right)^2+ O \left( \left( x - \dfrac{\pi}{2} \right)^3 \right) .
$$
Therefore, there exists a left neighborhood of $\pi/2$ such that:
\begin{equation}
\label{app_diff_Pi/2}
\varphi_{p, 2}(x) < 0
\mbox{ ,}\quad
\dfrac{\partial {\varphi_{p, 2}(x)}}{\partial x} > 0
\quad\mbox{and}\quad
\dfrac{\partial^2 {\varphi_{p, 2}(x)}}{\partial x^2} > 0  
\end{equation}
for $p\!\in\!(p_1, p_2)$.\\
By analyzing the monotonicity of the functions $\dfrac{\partial^3 {\varphi_{p, 2}(x)}}{\partial x^3}$, $\dfrac{\partial^2 {\varphi_{p, 2}(x)}}{\partial x^2}$, $\dfrac{\partial {\varphi_{p, 2}(x)}}{\partial x}$ and $\varphi_{p, 2}(x)$ for $p \!\in\! (p_1, p_2)$ on the interval $(0, \pi/2)$, in a similar manner as in the proof of Statement \ref{Statement_A}, based on (\ref{app_diff_3}), (\ref{app_diff_0}) and (\ref{app_diff_Pi/2}),
it can be concluded that each function  
$\varphi_{p, 2}(x)$, for $p \!\in\! (p_1, p_2)$, has exactly one minimum on the interval $(0, \pi/2)$. From $\displaystyle{\lim_{x \to 0+} \varphi_{p, 2}(x)} > 0$ and $\displaystyle{\lim_{x \to \pi/2-} \varphi_{p, 2}(x)} < 0$, for $p\!\in\!(p_1, p_2)$, the corresponding inequalities follow.

\medskip
\noindent
{\boldmath $(iii)$}
The claim follows directly from Statement \ref{Statement_A} and based on the stratification. Namely, for $q=2$, it holds that $A(q)=\dfrac{\pi-2}{\pi^{q+1}} = p_2$.

\medskip
\noindent
{\boldmath $(iv)$}
It has been proven within the proof {\boldmath $(ii)$}.

\medskip
\noindent
{\boldmath $(v)$}
Note that the infimum of the error $d(p) = \sup_{ x \in (0,\pi/2)}|\varphi_{p, 2}(x)|$, for $p\!\in\!(p_1, p_2)$, exists and is attained when:
\begin{equation}
\label{system_app1_1}
\left| \varphi_{p, 2} \left(0 + \right) \right| = \left| \varphi_{p, 2} \left( m^{(p)} \right) \right| .
\end{equation}
The equation (\ref{system_app1_1}) can be numerically solved using the Computer Algebra System \textit{Maple},  yielding in the value of the parameter $p = p_0$ being numerically determined as:
$$
p_0 = 0.036014\ldots \,,
$$
which determines the minimax approximant $\varphi_{p_0, 2}(x)$ of the family of functions $\varphi_{p, 2}(x)$.
\hfill $\Box$
\noindent

\medskip\noindent
In order to improve and extend Theorem \ref{Theorem_1_3}, we consider the family of functions $\varphi_{p, q}(x)$ for the case $q =3$. The family of functions $\varphi_{p, 3}(x)$ reduces to:
\begin{equation}
\label{family_T3}
\varphi_{p, 3}(x)
=
\dfrac{\sin x}{x} - \dfrac{2}{\pi} - p \left(\pi^3 - 8x^3 \right)
\end{equation}
and is decreasingly stratified with respect to the parameter $p \!\in\! R^{+}$ on the interval $(0, \pi/2)$, as proven in Lemma \ref{Lemma_1}.
For this family, the following statement holds:


\begin{statement}
\label{statement_T3}
Let$:$
$$
p_1
=
\dfrac{2}{3\pi^4}
=
0.0068439 \ldots
\quad\mbox{and}\quad
p_2
=
\dfrac{\pi - 2}{\pi^4}
=
0.011719 \ldots \,.
$$
Then, it holds$:$

\medskip
\noindent
{\boldmath $(i)$}
If $p\!\in\!(0,p_1]$, then$:$
$$
x \!\in\! \left(0, \dfrac{\pi}{2}\right)
\;\Longrightarrow\;
\dfrac{\sin x}{x}
>
\dfrac{2}{\pi} + p_1 \left(\pi^3 - 8x^3 \right)
\geq
\dfrac{2}{\pi} + p \left(\pi^3 - 8x^3 \right) .
$$

\medskip
\noindent
{\boldmath $(ii)$}
If $p\!\in\!(p_1, p_2)$, then the equality$:$
$$
\varphi_{p, 3}(x)
=
\dfrac{\sin x}{x} - \dfrac{2}{\pi} - p \left(\pi^3 - 8x^3 \right)
=
0
$$
has a unique solution $x_0^{(p)}$ and it holds$:$
$$
x\!\in\!\left( 0,x_{0}^{(p)} \right)
\;\Longrightarrow\;
\dfrac{\sin x}{x} > \dfrac{2}{\pi} + p \left(\pi^3 - 8x^3 \right)
$$
and
$$
x\!\in\!\left( x_{0}^{(p)}, \dfrac{\pi}{2}\right)
\;\Longrightarrow\;
\dfrac{\sin x}{x} < \dfrac{2}{\pi} + p \left(\pi^3 - 8x^3 \right) .
$$

\medskip
\noindent
{\boldmath $(iii)$}
If $p\!\in\![p_2,+\infty)$, then$:$
$$
x \!\in\! \left(0, \dfrac{\pi}{2}\right)
\;\Longrightarrow\;
\dfrac{\sin x}{x}
<
\dfrac{2}{\pi} + p_2 \left(\pi^3 - 8x^3 \right)
\leq
\dfrac{2}{\pi} + p \left(\pi^3 - 8x^3 \right)  .
$$

\medskip
\noindent
{\boldmath $(iv)$}
Each function from the family $\varphi_{p, 3}(x)$, for $p \!\in\! (p_1, p_2]$, has exactly one minimum at a point $m^{(p)} \!\in\! (0, \pi/2)$ on the interval $(0, \pi/2)$.

\medskip
\noindent
{\boldmath $(v)$}
The equality$:$
$$
\left| \varphi_{p, 3} \left(0 + \right) \right| = \left| \varphi_{p, 3} \left( m^{(p)} \right) \right|
$$
has the solution $p = p_0$, for the parameter $p \!\in\! (p_1, p_2)$, numerically determined as$:$
$$
p_0
=
0.010441 \ldots \,.
$$
For value$:$
$$
d_{0}
=
\left| \varphi_{p_0, 3} \left( 0 + \right) \right| = \left| \varphi_{p_0, 3} \left(m^{(p_0)} \right) \right|
=
0.039635 \ldots \,,
$$
it holds$:$
$$
d_0 = \inf\limits_{p \in (0,\infty)} \sup\limits_{x \in (0,\pi/2)}{\left| \varphi_{p, 3}(x) \right|} \, .
$$
Hence, the minimax approximant of the family of functions $\varphi_{p, 3}(x)$ is$:$
$$
\varphi_{p_0, 3}(x)
=
\dfrac{\sin x}{x} - \dfrac{2}{\pi} - p_0 \left(\pi^3 - 8x^3 \right) ,
$$
which determines the corresponding $($minimax$)$ approximation$:$
\begin{equation}
\label{minimax_approx_T3}
\dfrac{\sin x}{x}
\approx
\dfrac{2}{\pi} + 0.010441 \ldots \left(\pi^3 - 8x^3 \right) .
\end{equation}
\end{statement}
\noindent
{\bf Proof.} Analogously to the proof of Statement \ref{statement_T2}.
\hfill $\Box$

\medskip\noindent
In order to improve and extend Theorem \ref{Theorem_1_4}, we consider the family of functions $\varphi_{p, q}(x)$ for the case $q = 4$. The family of functions $\varphi_{p, 4}(x)$ reduces to:
\begin{equation}
\label{family_T4}
\varphi_{p, 4}(x)
=
\dfrac{\sin x}{x} - \dfrac{2}{\pi} - p\left(\pi^4 - 16x^4 \right)
\end{equation}
and is decreasingly stratified with respect to the parameter $p \!\in\! R^{+}$ on the interval $(0, \pi/2)$, as proven in Lemma \ref{Lemma_1}.
For this family, the following statement holds:


\begin{statement}
\label{statement_T4}
Let$:$
$$
p_1
=
\dfrac{1}{2\pi^5}
=
0.0016338 \ldots
\quad\mbox{and}\quad
p_2
=
\dfrac{\pi - 2}{\pi^5}
=
0.0037304 \ldots \,.
$$
Then, it holds$:$

\medskip
\noindent
{\boldmath $(i)$}
If $p\!\in\!(0,p_1]$, then$:$
$$
x \!\in\! \left(0, \dfrac{\pi}{2}\right)
\;\Longrightarrow\;
\dfrac{\sin x}{x}
>
\dfrac{2}{\pi} + p_1 \left(\pi^4 - 16x^4 \right)
\geq
\dfrac{2}{\pi} + p \left(\pi^4 - 16x^4 \right) .
$$

\medskip
\noindent
{\boldmath $(ii)$}
If $p\!\in\!(p_1, p_2)$, then the equality$:$
$$
\varphi_{p, 4}(x)
=
\dfrac{\sin x}{x} - \dfrac{2}{\pi} - p \left(\pi^4 - 16x^4 \right)
=
0
$$
has a unique solution $x_0^{(p)}$ and it holds$:$
$$
x\!\in\!\left( 0,x_{0}^{(p)} \right)
\;\Longrightarrow\;
\dfrac{\sin x}{x} > \dfrac{2}{\pi} + p \left(\pi^4 - 16x^4 \right)
$$
and
$$
x\!\in\!\left( x_{0}^{(p)}, \dfrac{\pi}{2}\right)
\;\Longrightarrow\;
\dfrac{\sin x}{x} < \dfrac{2}{\pi} + p \left(\pi^4 - 16x^4 \right) .
$$

\medskip
\noindent
{\boldmath $(iii)$}
If $p\!\in\![p_2,+\infty)$, then$:$
$$
x \!\in\! \left(0, \dfrac{\pi}{2}\right)
\;\Longrightarrow\;
\dfrac{\sin x}{x}
<
\dfrac{2}{\pi} + p_2 \left(\pi^4 - 16x^4 \right)
\leq
\dfrac{2}{\pi} + p \left(\pi^4 - 16x^4 \right) .
$$

\break

\medskip
\noindent
{\boldmath $(iv)$}
Each function from the family $\varphi_{p, 4}(x)$, for $p \!\in\! (p_1, p_2]$, has exactly one minimum at a point $m^{(p)} \!\in\! (0, \pi/2)$ on the interval $(0, \pi/2)$.

\medskip
\noindent
{\boldmath $(v)$}
The equality$:$
$$
\left| \varphi_{p, 4} \left(0 + \right) \right| = \left| \varphi_{p, 4} \left( m^{(p)} \right) \right|
$$
has the solution $p = p_0$, for the parameter $p \!\in\! (p_1, p_2)$, numerically determined as$:$
$$
p_0
=
0.0031146 \ldots \,.
$$
For value$:$
$$
d_{0}
=
\left| \varphi_{p_0, 4} \left( 0 + \right) \right| = \left| \varphi_{p_0, 4} \left(m^{(p_0)} \right) \right|
=
0.059981 \ldots \,,
$$
it holds$:$
$$
d_0 = \inf\limits_{p \in (0,\infty)} \sup\limits_{x \in (0,\pi/2)}{\left| \varphi_{p, 4}(x) \right|} \, .
$$
Hence, the minimax approximant of the family of functions $\varphi_{p, 4}(x)$ is$:$
$$
\varphi_{p_0, 4}(x)
=
\dfrac{\sin x}{x} - \dfrac{2}{\pi} - p_0 \left(\pi^4 - 16x^4 \right) ,
$$
which determines the corresponding $($minimax$)$ approximation$:$
\begin{equation}
\label{minimax_approx_T4}
\dfrac{\sin x}{x}
\approx
\dfrac{2}{\pi} + 0.0031146 \left(\pi^4 - 16x^4 \right) .
\end{equation}
\end{statement}
\noindent
{\bf Proof.} Analogously to the proof of Statement \ref{statement_T2}.
\hfill $\Box$

\bigskip\noindent
Figure 3 illustrates the stratified families of functions $\varphi_{p, 1}(x)$, $\varphi_{p, 2}(x)$, $\varphi_{p, 3}(x)$ and $\varphi_{p, 4}(x)$ respectively, see $(\ref{family_T1})$, $(\ref{family_T2})$, $(\ref{family_T3})$ and $(\ref{family_T4})$.
For each family, cases for all values of the parameter $p \!\in\! R^+$ are shown. Particularly, cases with constants obtained in Statement \ref{statement_T1}, \ref{statement_T2}, \ref{statement_T3} and \ref{statement_T4}, some of which are also obtained in Theorems  \ref{Theorem_1_1}, \ref{Theorem_1_2}, \ref{Theorem_1_3} and \ref{Theorem_1_4}, are singled out.

\break

\medskip\noindent
\begin{figure}[htbp]
\centering
	\begin{subfigure}[b]{0.45\textwidth}
            \centering
            \includegraphics[width=\textwidth]{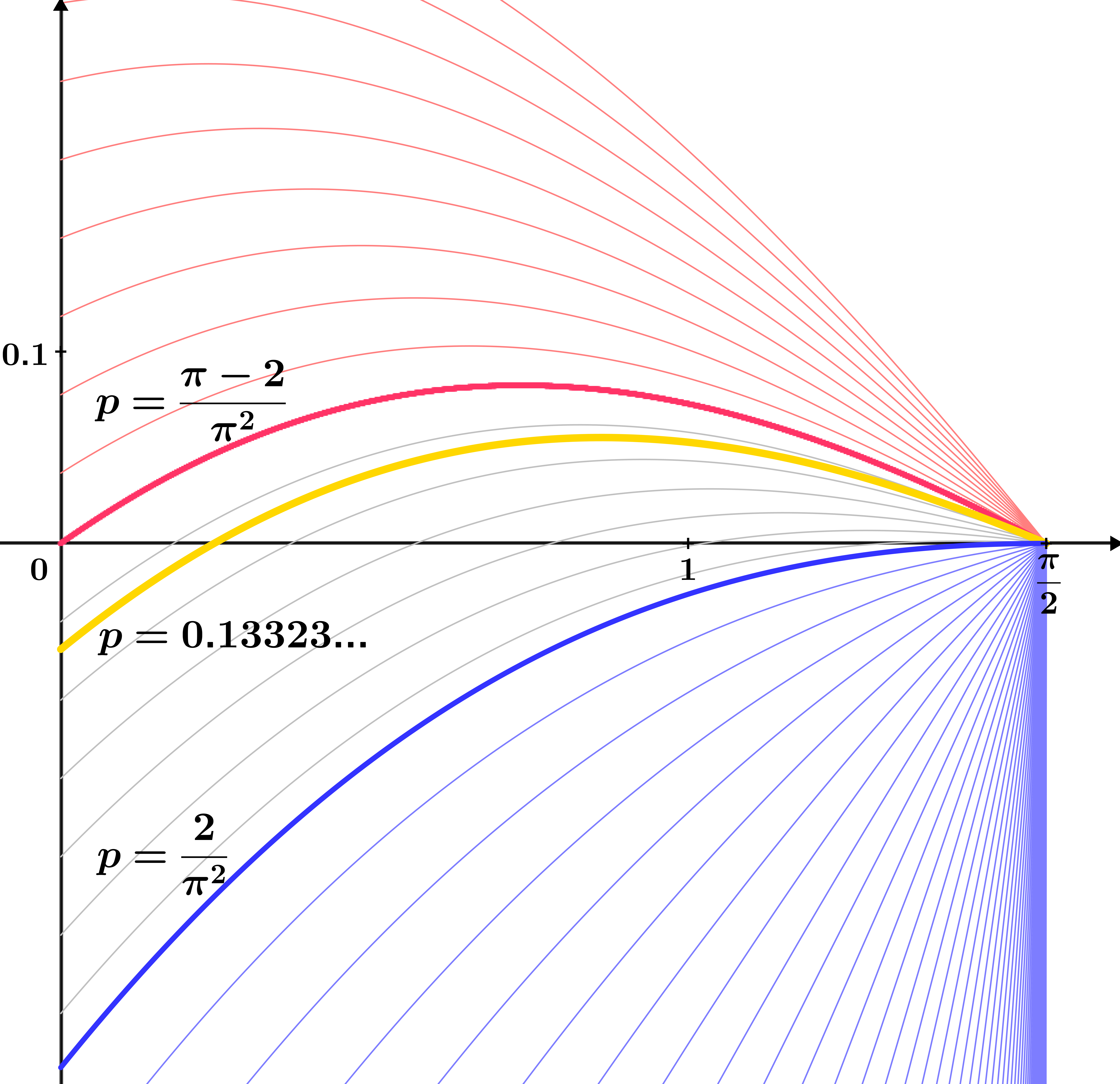}
            \caption[Network2]%
            {{\small $\varphi_{p, 1}(x)$, see $(\ref{family_T1})$}}    
        \end{subfigure}
        \hfill
        \begin{subfigure}[b]{0.45\textwidth}  
            \centering 
            \includegraphics[width=\textwidth]{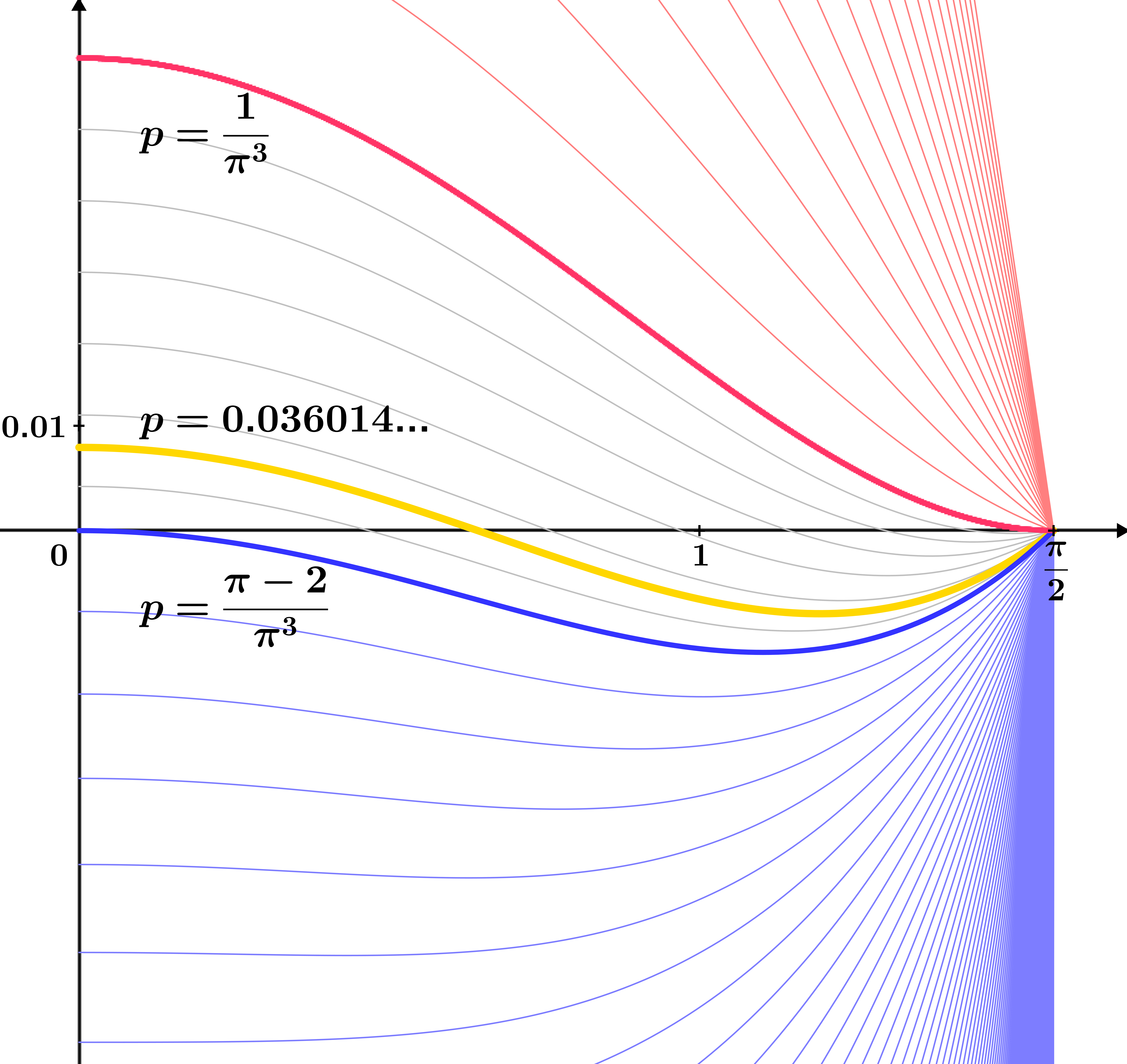}
            \caption[]%
            {{\small $\varphi_{p, 2}(x)$, see $(\ref{family_T2})$}}
        \end{subfigure}
        \vskip\baselineskip
        \begin{subfigure}[b]{0.45\textwidth}   
            \centering 
            \includegraphics[width=\textwidth]{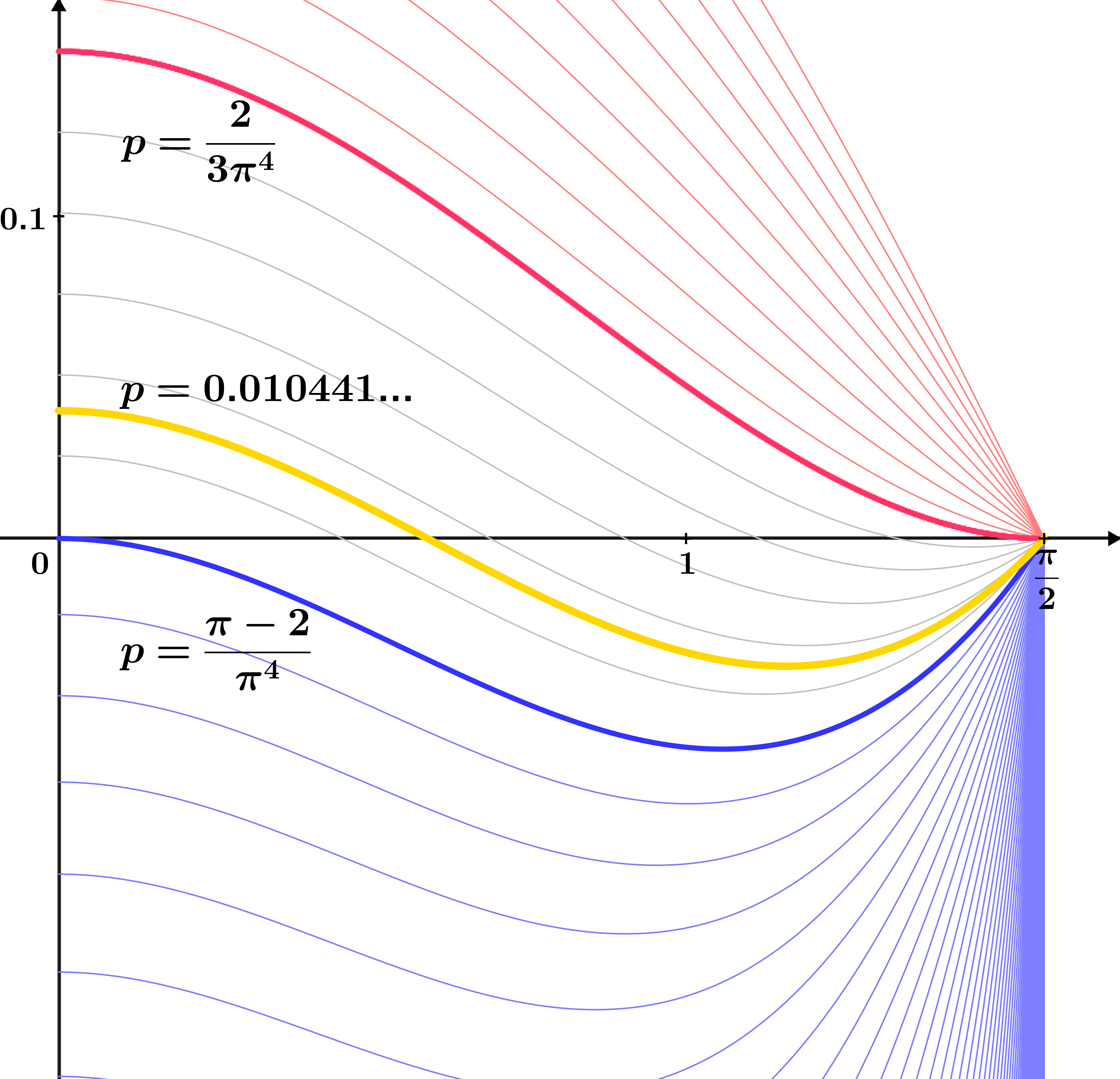}
            \caption[]%
            {{\small $\varphi_{p, 3}(x)$, see $(\ref{family_T3})$}}
        \end{subfigure}
        \hfill
        \begin{subfigure}[b]{0.45\textwidth}   
            \centering 
            \includegraphics[width=\textwidth]{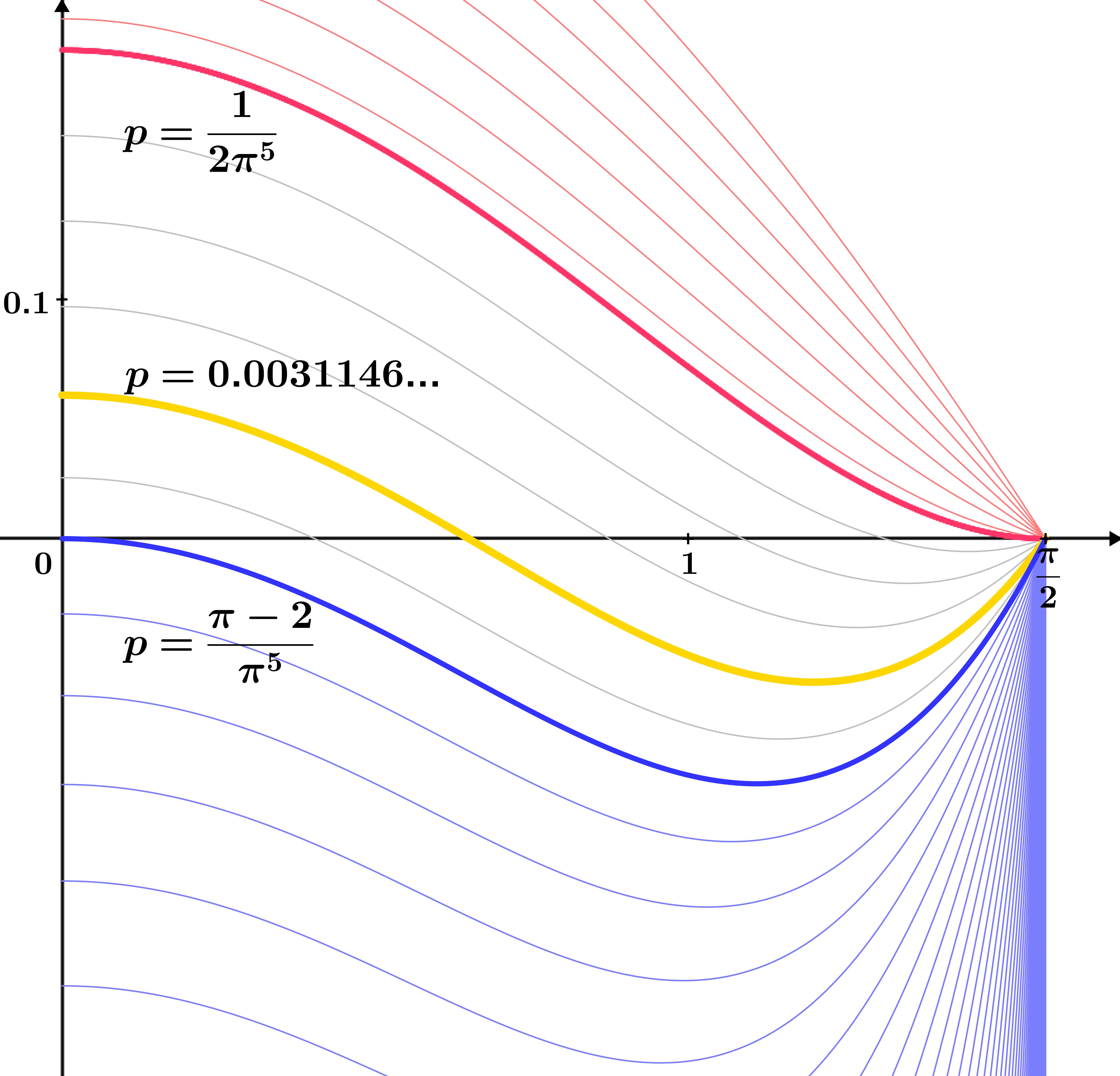}
            \caption[]%
            {{\small $\varphi_{p, 4}(x)$, see $(\ref{family_T4})$}}   
        \end{subfigure}
\caption{Stratified families of functions: ({a}) $\varphi_{p, 1}(x)$, ({b}) $\varphi_{p, 2}(x)$, ({c}) $\varphi_{p, 3}(x)$, ({d}) $\varphi_{p, 4}(x)$ \label{Jordan-Family T1234}}
\end{figure}

\vspace{-1mm}

\subsection{Approximations of the ${\mbox{\boldmath sinc}}$ function}

\medskip\noindent
In this subsection, we provide some approximations of the $\mbox{\rm sinc}$ function and analyze the maximum approximation errors.
The previously obtained upper and lower bounds of the $\mbox{\rm sinc}$  function can be used to derive some approximations of this function. Further, more optimal approximations can be obtained through the corresponding minimax approximants. 

\medskip\noindent
In Table \ref{table_upper_bounds}, we present some upper bounds of the $\mbox{\rm sinc}$ function derived from Theorems \ref{Theorem_1_1}, \ref{Theorem_1_2}, \ref{Theorem_1_3} and \ref{Theorem_1_4}, that is, Statements \ref{statement_T1}, \ref{statement_T2}, \ref{statement_T3} and \ref{statement_T4} and Statements \ref{Statement_A} and \ref{Statement_B}. It is noteworthy that the upper bound from Theorem \ref{Theorem_1_2} (the best upper bound from Statement \ref{statement_T2}) is identical to the best upper bound from Statement \ref{Statement_A}.

\begin{table}[hbt!]
\vspace{-1mm}
\hspace*{-0.55cm}
\begin{tabular} {||c|c||}
 \hline
 Upper bound                                                 & Maximum deviation\\ [0.2ex] 
 of the $\mbox{\rm sinc}\,x$ function & from the $\mbox{\rm sinc}\,x$ function\\ [0.2ex] 
 on the interval $(0, \pi/2)$                                                & on the interval $(0, \pi/2)$\\ [0.01ex]
 \hline\hline
\rule{0pt}{4.7ex} $\dfrac{\sin x}{x} < \dfrac{2}{\pi} + \dfrac{2}{\pi^2}\left(\pi - 2x \right)$ & $\dfrac{\pi-4}{\pi} = 0.27323 \ldots$ \\ [2.5ex]
 \hline
\rule{0pt}{4.7ex} $\dfrac{\sin x}{x} < \dfrac{2}{\pi} + \dfrac{\pi - 2}{\pi^3}\left(\pi^2 - 4x^2 \right)$ & $0.011612 \ldots$ \\ [2.5ex] 
 \hline
\rule{0pt}{4.7ex} $\dfrac{\sin x}{x} < \dfrac{2}{\pi} + \dfrac{\pi - 2}{\pi^4}\left(\pi^3 - 8x^3 \right)$ & $0.065358 \ldots$ \\ [2.5ex] 
 \hline
\rule{0pt}{4.7ex} $\dfrac{\sin x}{x} < \dfrac{2}{\pi} + \dfrac{\pi - 2}{\pi^5}\left(\pi^4 - 16x^4 \right)$ & $0.10245 \ldots$ \\ [2.5ex] 
 \hline
\rule{0pt}{4.7ex} $\dfrac{\sin x}{x} < \dfrac{2}{\pi} + \dfrac{2}{  \left({\frac{\pi^2}{4}-1}\right) \pi^{\frac{\pi^2}{4}}} \left( \pi^{\frac{\pi^2}{4}-1} - \left(2x\right)^{\frac{\pi^2}{4}-1} \right)$ & $\dfrac{-\pi^2+2\pi+4}{\pi^2-4} = 0.070461 \ldots$ \\ [2.5ex] 
 \hline
\end{tabular}
\centering
\caption{\label{table_upper_bounds} Upper bounds of the $\mbox{\rm sinc}\,x$ function on the interval $(0, \pi/2)$}
\end{table}

\vspace{-1mm}

\noindent
In Table \ref{table_lower_bounds}, we present some lower bounds of the $\mbox{\rm sinc}$ function derived from Theorems \ref{Theorem_1_1}, \ref{Theorem_1_2}, \ref{Theorem_1_3} and \ref{Theorem_1_4}, that is, Statements \ref{statement_T1}, \ref{statement_T2}, \ref{statement_T3} and \ref{statement_T4} and Statements \ref{Statement_A} and \ref{Statement_B}. It is noteworthy that the best lower bound from Statement \ref{Statement_A} is identical to the best lower bound from Statement \ref{Statement_B}.

\begin{table}[hbt!]
\begin{tabular} {||c|c||}
 \hline
 Lower bound                                                 & Maximum deviation\\ [0.2ex] 
 of the $\mbox{\rm sinc}\,x$ function & from the $\mbox{\rm sinc}\,x$ function\\ [0.2ex] 
 on the interval $(0, \pi/2)$                                                & on the interval $(0, \pi/2)$\\ [0.01ex]
 \hline\hline
\rule{0pt}{4.7ex} $\dfrac{2}{\pi} + \dfrac{\pi-2}{\pi^2} \left(\pi - 2x \right) <\dfrac{\sin x}{x}$ & $0.082395 \ldots$ \\ [2.5ex] 
 \hline
\rule{0pt}{4.7ex} $\dfrac{2}{\pi} + \dfrac{1}{\pi^3} \left(\pi^2 - 4x^2 \right) <\dfrac{\sin x}{x}$ & $\dfrac{\pi-3}{\pi} = 0.045070 \ldots$ \\ [2.5ex] 
 \hline
\rule{0pt}{4.7ex} $\dfrac{2}{\pi} + \dfrac{2}{3\pi^4} \left(\pi^3 - 8x^3 \right) <\dfrac{\sin x}{x}$ & $\dfrac{3\pi-8}{3\pi} = 0.15117 \ldots$ \\ [2.5ex] 
 \hline
\rule{0pt}{4.7ex} $\dfrac{2}{\pi} + \dfrac{1}{2\pi^5} \left(\pi^4 - 16x^4 \right) <\dfrac{\sin x}{x}$ & $\dfrac{2\pi-5}{2\pi} = 0.20422 \ldots$ \\ [2.5ex] 
 \hline
\rule{0pt}{4.7ex} $\dfrac{2}{\pi} + \dfrac{\pi-2}{\pi^{\frac{2}{\pi-2}+1}} \left( \pi^{\frac{2}{\pi-2}} - \left(2x\right)^{\frac{2}{\pi-2}} \right) < \dfrac{\sin x}{x}$ & $0.0085153 \ldots$ \\ [2.5ex] 
 \hline
\end{tabular}
\centering
\caption{\label{table_lower_bounds} Lower bounds of the $\mbox{\rm sinc}\,x$ function on the interval $(0, \pi/2)$}
\end{table}

\break

\medskip\noindent
In Table  \ref{table_minimax_approximations}, we present some minimax approximations of the $\mbox{\rm sinc}$ function derived from  the minimax approximants of the families $\varphi_{p, 1}(x)$, $\varphi_{p, 2}(x)$, $\varphi_{p, 3}(x)$, $\varphi_{p, 4}(x)$, $\varphi_{A(q), q}(x)$ and $\varphi_{B(q), q}(x)$ respectively. These families are considered in Statements \ref{statement_T1}, \ref{statement_T2}, \ref{statement_T3} and \ref{statement_T4} with the aim of improving Theorems \ref{Theorem_1_1}, \ref{Theorem_1_2}, \ref{Theorem_1_3} and \ref{Theorem_1_4}, respectively, and in Statements \ref{Statement_A} and \ref{Statement_B}.

\begin{table}[hbt!]
\hspace*{-0.5cm}
\begin{tabular} {||c|c||}
 \hline
 Minimax approximation                                                 & Maximum deviation\\ [0.2ex] 
 of the $\mbox{\rm sinc}\,x$ function & from the $\mbox{\rm sinc}\,x$ function\\ [0.2ex] 
 on the interval $(0, \pi/2)$                                               & on the interval $(0, \pi/2)$\\ [0.01ex]
 \hline\hline
\rule{0pt}{4.7ex} $\dfrac{\sin x}{x} \approx \dfrac{2}{\pi} + 0.13323\ldots \left(\pi - 2x \right)$ & $0.055187 \ldots$ \\ [2.5ex] 
 \hline
\rule{0pt}{4.7ex} $\dfrac{\sin x}{x} \approx \dfrac{2}{\pi} + 0.036014\ldots \left(\pi^2 - 4x^2 \right)$ & $0.0079283 \ldots$ \\ [2.5ex] 
 \hline
\rule{0pt}{4.7ex} $\dfrac{\sin x}{x} \approx \dfrac{2}{\pi} + 0.010441\ldots \left(\pi^3 - 8x^3 \right)$ & $0.039635 \ldots$ \\ [2.5ex] 
 \hline
\rule{0pt}{4.7ex} $\dfrac{\sin x}{x} \approx \dfrac{2}{\pi} + 0.0031146\ldots \left(\pi^4 - 16x^4 \right)$ & $0.059981 \ldots$ \\ [2.5ex] 
 \hline
\rule{0pt}{4.7ex} $\dfrac{\sin x}{x} \approx \dfrac{2}{\pi} + 0.043803 \ldots \left( \pi^{1.84823 \ldots} - \left(2x\right)^{1.84823 \ldots} \right)$ & $0.0026604 \ldots$ \\ [2.5ex] 
 \hline
\rule{0pt}{4.7ex} $\dfrac{\sin x}{x} \approx \dfrac{2}{\pi} + 0.051415 \ldots \left( \pi^{1.72287  \ldots} - \left(2x\right)^{1.72287  \ldots} \right)$ & $0.0061296 \ldots$ \\ [2.5ex] 
 \hline
\end{tabular}
\centering
\caption{\label{table_minimax_approximations} Minimax approximations of the $\mbox{\rm sinc}\,x$ function on the interval $(0, \pi/2)$}
\end{table}

\section{Conclusion}

\medskip\noindent
In this paper, two double Jordan-type inequalities have been obtained, encompassing the inequalities established in the papers  \cite{Qi_Guo_1993}-\cite{Li_Li_2008}}. These inequalities were explored in the context of stratified families of functions, a concept introduced in recent research \cite{Malesevic_Mihailovic_2021}. The introduction of stratified families of functions enables the derivation of known results for specific parameter choices,  including the analysis of parameter values previously unknown in the Theory of Analytic Inequalities. Furthermore, we identify parameter values within each examined family of functions for which the function, as a member of that family, exhibits some optimal properties (minimax approximant). Based on these minimax approximants and functions representing the upper and lower bounds of the $\mbox{\rm sinc}$ function, we provided some approximations of the $\mbox{\rm sinc}$ function. Additionally, we analyzed the errors associated with all mentioned approximations.

\medskip\noindent
It is crucial to emphasize that the minimax approximant of the stratified family of functions is the function for which the minimal error in approximations is obtained within the given family of functions. Therefore, identifying these parameter values is significant in the Approximation Theory.

\medskip\noindent
By considering the stratified family of functions individually with respect to two parameters, we were able to analyze Jordan-type inequalities in a unified manner, resulting in both previously established and novel findings. Future research endeavors will focus on extending this approach even further.

\bigskip
\noindent
{\bf Acknowledgments.} The work of the authors was financially supported by the Ministry of Science, Technological Development and Innovation of the Republic of Serbia under contract numbers: 451-03-66/2024-03/200103 (for the first author) and 451-03-65/2024-03/200103 (for the second author).
 \\[0.5 ex]

\vspace{1cc}
{\small
\noindent
\textbf{Milo\v s Mi\' covi\' c} \\
University of Belgrade, \\
School of Electrical Engineering, \\
Bulevar kralja Aleksandra 73, \\
11000 Belgrade, Serbia, \\
E-mail: {\it milos.micovic@etf.bg.ac.rs}

\vspace{1cc}
{\small
\noindent
\textbf{Branko Male\v {s}evi\' c} \\
University of Belgrade, \\
School of Electrical Engineering, \\
Bulevar kralja Aleksandra 73, \\
11000 Belgrade, Serbia, \\
E-mail: {\it branko.malesevic@etf.bg.ac.rs}


\begin{thebibliography}{5}

\bibitem{Qi_Guo_1993}
F. Qi and B.-N. Guo, {\em  On generalizations of Jordan's inequality} (in Chinese), Coal Higher Education, supplement, 1993,
32-33

\bibitem{Qi_1996}
{\sc  F. Qi}, {\em  Extensions and sharpenings of Jordan's and Kober's inequality} (in Chinese),
Journal of Mathematics for Technology 1996, 12:4, 98-102

\bibitem{Deng_1995}
{\sc  K. Deng}, {\em  The noted Jordan's inequality and its extensions} (in Chinese),
Journal of Xiangtan Mining Institute 1995, 10:4, 60-63

\bibitem{Jiang_Yun_2006}
{\sc  W. D. Jiang}, {\sc H. Yun}, {\em Sharpening of Jordan's inequality and its applications},
Journal of Inequalities in Pure and Applied Mathematics 2006, 7:3, 1-4

\bibitem{Li_Li_2008}
{\sc  J.-L. Li}, {\sc Y.-L. Li}, {\em On the Strengthened Jordan's Inequality},
Journal of Inequalities and Applications 2008, 2007, 1-8

\bibitem{Malesevic_Mihailovic_2021}
{\sc B. Male\v sevi\' c}, {\sc B. Mihailovi\' c}, {\em A minimax approximant in the theory of analytic inequalities},
Applicable Analysis and Discrete Mathematics 2021, 15:2, 486-509

\bibitem{Mitrinovic_1970}
{\sc D. Mitrinovi\' c}, {\em Analytic Inequalities}, Springler-Verlag 1970

\bibitem{Qi_Niu_Guo_1996}
{\sc  F. Qi}, {\sc  D.-W. Niu}, {\sc  B.-N. Guo}, {\em  Refinements, Generalizations, and Applications of
Jordan's Inequality and Related Problems}, 
Journal of Inequalities and Applications 2009, 2009, 1-52

\bibitem{Ozban_2006}
{\sc A. Y.\"{O}zban}, {\em A new refined form of Jordan's inequality and its applications}, Applied Mathematics Letters 2006, 19:2, 155-160

\bibitem{Li_2006}
{\sc  J.-L. Li}, {\em An identity related to Jordan's inequality},
International Journal of Mathematics and Mathematical Sciences 2006, 2006, 1-6

\bibitem{Zhu_2006}
{\sc L. Zhu}, {\em Sharpening Jordan's inequality and the Yang Le inequality},
Applied Mathematics Letters 2006, 19:3, 240-243

\bibitem{Zhu2_2006}
{\sc L. Zhu}, {\em Sharpening Jordan's inequality and the Yang Le inequality, II},
Applied Mathematics Letters 2006, 19:9, 990-994

\bibitem{Zhu_2008}
{\sc L. Zhu}, {\em A general refinement of Jordan-type inequality},
Computers \& Mathematics with Applications 2008, 55:11, 2498-2505

\bibitem{Niu_Huo_Cao_Qi_2008}
{\sc D.-W. Niu}, {\sc Z.-H. Huo}, {\sc J. Cao}, {\sc F. Qi},
{\em A general refinement of Jordan's inequality and a refinement of L. Yang's inequality},
Integral Transforms and Special Functions 2008, 19:3, 157-164

\bibitem{Chen_Debnath_2012}
{\sc C.-P. Chen}, {\sc L. Debnath},
{\em Sharpness and generalization of Jordan's inequality and its application},
Applied Mathematics Letters 2012, 25:3, 594-599

\bibitem{Barbu_Elias_2014}
{\sc C. Barbu}, {\sc L.-I. Pi\c{s}coran}, {\em Jordan type inequalities using monotony of functions}, Journal of Mathematical Inequalities 2014, 8:1, 83-89

\bibitem{Aharonov_Piscoran_2014}
{\sc D. Aharonov}, {\sc U. Elias}, {\em More Jordan type inequalities}, 
Mathematical Inequalities \& Applications 2014, 17:4, 1563-1577 

\bibitem{Alzer_Kwong_2018}
{\sc H. Alzer}, {\sc M. K. Kwong},
{\em On Jordan's inequality},
Periodica Mathematica Hungarica 2018, 77:2, 191-200

\bibitem{Zhang_Ma_2018}
{\sc L. Zhang}, {\sc X. Ma},
{\em New Refinements and Improvements of Jordan's Inequality},
Mathematics 2018, 6:12, 1-8

\bibitem{Zhang_Ma_2019}
{\sc L. Zhang}, {\sc X. Ma},
{\em New Polynomial Bounds for Jordan's and Kober's Inequalities Based on the Interpolation and Approximation Method},
Mathematics 2019, 7:8, 1-9

\bibitem{Zhang_Chen_2019}
{\sc B. Zhang}, {\sc C.-P. Chen}, {\em Sharpness and generalization of Jordan, Becker-Stark and Papenfuss inequalities with an application}, Journal of Mathematical Inequalities 2019, 13:4, 1209-1234

\bibitem{Haque_2020}
{\sc N. Haque},
{\em A Short Calculus Proof of Jordan's Inequality},
Cambridge Open Engage (Preprint) 2020, 1-1

\bibitem{Popa_2020}
{\sc E. C. Popa},
{\em A note on Jordan's inequality},
General Mathematics 2020, 28:2, 97-102

\bibitem{Malesevic_Micovic_2023}
{\sc B. Male\v sevi\' c}, {\sc M. Mi\' covi\' c}, {\em Exponential Polynomials and Stratification in the Theory of Analytic Inequalities}, 
Journal of Science and Arts 2023, 23:3, 659-670

\bibitem{Malesevic_Jovanovic_2024}
{\sc B. Male\v sevi\' c}, {\sc D. Jovanovi\' c}, {\em Frame's Types of Inequalities and Stratification}, CUBO, A Mathematical Journal 2024, 26:1, 1-19

\bibitem{Chen_Ge_2021}
{\sc S. Chen}, {\sc X. Ge}, {\em A solution to an open problem for Wilker-type inequalities}, 
Journal of Mathematical Inequalities 2021, 15:1, 59-65

\bibitem{Malesevic_Mihailovic_NenezicJovic_Milinkovic_2022}
{\sc B. Male\v sevi\' c}, {\sc B. Mihailovi\' c}, {\sc M. Nenezi\' c Jovi\' c}, {\sc L. Milinkovi\' c},
{\em Some minimax approximants of D'Aurizio trigonometric inequalities}, 
HAL (Preprint) 2022, 1-9

\bibitem{Chen_Mortici_2023}
{\sc C.-P. Chen}, {\sc C. Mortici}, {\em The relationship between Huygens' and Wilker's inequalities and further remarks},
Applicable Analysis and Discrete Mathematics 2023, 17:1, 92-100

\bibitem{Sandor_2016}
{\sc J. S\' andor}, {\em On D'Aurizio's trigonometric inequality}, 
Journal of Mathematical Inequalities 2016, 10:3, 885-888 

\bibitem{Sandor_2017}
{\sc J. S\' andor}, {\em Extensions of D'Aurizio's trigonometric inequality}, 
Notes on Number Theory and Discrete Mathematics 2017, 23:2, 81-83 

\bibitem{Hung_Li_2016}
{\sc L.-C. Hung}, {\sc P.-Y. Li}, {\em On generalization of D'Aurizio-S\' andor inequalities involving a parameter}, 
Journal of Mathematical Inequalities 2018, 12:3, 853-860 

\bibitem{Malesevic_Makragic_2016}
{\sc B. Male\v sevi\' c}, {\sc M. Makragi\' c}, {\em A method for proving some inequalities on mixed trigonometric polynomial functions}, Journal of Mathematical Inequalities 2016, 10:3, 849-876

\bibitem{Pinelis_2002}
{\sc I. Pinelis}, {\em L'Hospital type rules for monotonicity, with applications}, Journal of Inequalities in Pure and
Applied Mathematics 2002, 3:1, 1-5

\bibitem{Estrada_Pavlovic_2017}
{\sc R. Estrada}, {\sc M. Pavlovi\' c}, {\em L'H\^ opital's monotone rule, Gromov's theorem, and operations that preserve the monotonicity of quotients}, Publications de l'Institut Mathematique 2017, 101:115, 11-24

\bibitem{Cutland_1980}
{\sc N. Cutland}, {\em Computalibity: an introduction to recursive funtion theory},
Cambridge University Press 1980

\end{thebibliography}
\end{document}